\input amstex.tex
\documentstyle{amsppt}
\pageno=1
\hyphenpenalty=200
\nologo

\overfullrule=0pt
\hyphenation{semi-simple}

\TagsOnRight
\define\ga{\alpha}
\define\gb{\beta}
\define\gc{\gamma}
\define\gC{\Gamma}

\define\gD{\Delta}

\define\gth{\theta}

\define\gk{\kappa}
\define\gl{\lambda}

\define\gm{\mu}

\define\gr{\rho}
\define\gs{\sigma}

\define\gt{\tau}

\define\gw{\omega}
\define \gW{\Omega}

\title Cohomologie Automorphe et Compactifications partielles  de 
certaines Vari\'et\'es de Griffiths--Schmid. \endtitle

\author Henri Carayol  \endauthor

\address {I.R.M.A., 7 Rue Ren\'e Descartes, 67084 Strasbourg
Cedex, France}\endaddress
\email carayol\@math.u-strasbg.fr\endemail

\leftheadtext{Cohomologie Automorphe}
\rightheadtext{et Compactifications partielles \dots}



\abstract
{We study some automorphic cohomology classes of degree one on the Griffiths-Schmid varieties
attached  to some unitary groups in 3 variables.  Using   partial compactifications  of those varieties,  constructed by K. Kato and S. Usui,  
 we define for such a cohomology class  some analogues   of 
Fourier-Shimura coefficients,  
 which are cohomology classes on certain elliptic curves.  We show  that a large space of such  automorphic classes  can
 be generated by those with
rational "coefficients".   More precisely, we consider those cohomology classes that come
from Picard modular forms, via some Penrose-like transform studied in a previous article :  we prove that
the coefficients of the classes thus obtained can be computed from  the coefficients of the Picard
form by a similar transform defined at the level of the elliptic curve.
 }
\endabstract

\endtopmatter

\document

Le pr\'esent travail prolonge les articles  \cite{Ca1} et \cite{Ca2}, dans
lesquels on \'etudiait les groupes de cohomologie de certains faisceaux 
coh\'erents sur les vari\'et\'es de Griffiths--Schmid associ\'ees \`a certaines 
formes de groupes unitaires en trois variables. Rappelons sommairement la d\'efinition des
vari\'et\'es en  question :  elles sont  des quotients de domaines de p\'eriodes (''period 
matrix domains''), attach\'es \`a certains groupes r\'eductifs, par des groupes 
arithm\'etiques, ou,
si l'on pr\'ef\`ere, des  version ad\'eliques de tels quotients. En tant que telles,
elles g\'en\'eralisent la notion de vari\'et\'e de Shimura, bien qu'elles ne soient pas
en g\'en\'eral, contrairement \`a ces derni\`eres, des vari\'et\'es alg\'ebriques.
Elles avaient
\'et\'e  consid\'er\'ees et \'etud\'ees par Griffiths et Schmid (cf. \cite{G-S}) parce
qu'elles param\`etrent des structures de Hodge munies de polarisations et de donn\'ees 
additionnelles.
\medskip
Comme dans le cas plus classique des vari\'et\'es de Shimura, la d\'ecomposition 
des  groupes de
cohomologie de faisceaux coh\'erents se d\'ecrit en termes de formes automorphes, une
diff\'erence essentielle avec ce cas habituel consistant en ceci : que le $H^0$ est 
souvent nul ou trivial, autrement dit qu'il n'y a pas ou peu de formes automorphes au
sens usuel (d'o\`u le fait que ces vari\'et\'es sont la plupart du temps non-alg\'ebriques),
mais plut\^ot des classes de cohomologie automorphe de degr\'e $\geq 1$.  Les
articles  \cite{Ca1} et \cite{Ca2} tirent leur origine de l'observation que ces groupes de cohomologie,
dans le cas des vari\'et\'es de Griffiths--Schmid, peuvent contenir un peu plus
de formes automorphes que dans le
cas des vari\'et\'es de Shimura; en particulier, du moins dans le cas de formes de groupes
unitaires
\`a trois variables, mais sans doute aussi dans des cas beaucoup plus g\'en\'eraux, on y
trouve des classes associ\'ees \`a des formes du type de Maass qui ne peuvent pas intervenir
dans le cadre des vari\'et\'es de Shimura. On aimerait bien s\^ur beaucoup
pouvoir se servir de l'apparition de ces formes dans la cohomologie de ces
vari\'et\'es afin de prouver des propri\'et\'es arithm\'etiques (rationalit\'e,
congruences, construction de repr\'esentations galoisiennes), mais on se heurte 
aussit\^ot dans ce projet \`a l'obstacle essentiel que constitue leur non-alg\'ebricit\'e.
\medskip

L'objet du pr\'esent article est d'aborder la  question suivante : comment peut-on d\'efinir ce
qu'est une classe de cohomologie rationnelle, en d\'epit du fait que la vari\'et\'e \'etudi\'ee
n'est pas munie d'une structure arithm\'etique ? Si l'on pense au cas plus habituel des
vari\'et\'es de Shimura (en oubliant provisoirement qu'elles sont alg\'ebriques et d\'efinies sur
 un corps de nombres)
 et si l'on veut d\'efinir \'el\'emen\-tairement 
 pour ces derni\`eres ce qu'est la rationalit\'e d'une forme
automorphe, on voit que cela peut se faire essentiellement de deux fa\c cons : 
la premi\`ere consiste \`a regarder
les valeurs prises aux points sp\'eciaux et la seconde \`a consid\`erer 
les coefficients des d\'eveloppements
 de Fourier aux pointes ou composantes fronti\`eres. La premi\`ere fa\c con peut se g\'en\'eraliser sous la
forme suivante aux classes de cohomologie port\'ees par les vari\'et\'es de 
Griffiths-Schmid : on int\`egre ladite classe sur des cycles convenables, ce qui est
\`a la base d'une approche souvent utilis\'ee (cf. par exemple \cite{WW}) quoique 
pour d'autres motifs que des questions de rationalit\'e.  Quant \`a la seconde fa\c con
de proc\'eder, sa g\'en\'eralisation au cas pr\'esent n\'ecessite de disposer d'un analogue 
des composantes fronti\`eres, c'est-\`a-dire d'une sorte de compactification des 
vari\'et\'es de Griffiths-Schmid. Or une telle compactification partielle vient d'\^etre 
construite par K. Kato et S. Usui \cite{K-U}. Le  but du pr\'esent article  est d'explorer
ce que l'on peut faire \`a l'aide de ces compactifications.
\medskip

Nous nous pla\c cons ici dans le cas
particulier d'un groupe unitaire en trois va\-riables ;  c'est le m\^eme cadre que celui o\`u
nous nous
\'etions plac\'es dans les deux articles \cite{Ca1} \cite{Ca2},  \`a ceci pr\`es que nous prenions
alors pour simplifier une forme anisotrope du groupe, ce que nous ne faisons plus maintenant
car bien s\^ur nous voulons travailler sur une vari\'et\'e non-compacte. Notre vari\'et\'e complexe est
de dimension 3,  et nous d\'ecrivons dans ce cas la "compactification" (en fait tr\`es
partielle) que nous donne la construction de  Kato et Usui : elle revient \`a ajouter comme
composantes fronti\`eres des courbes elliptiques \`a multiplication complexe.
 D'autre part nous appliquons  une construction analogue \`a celle 
de \cite{Ca2} qui donne, \`a partir de formes de Picard classiques, les 
classes de cohomologie  - de degr\'e 1 - qui leur correspondent sur notre vari\'et\'e.
 Nous regardons ensuite la restriction au voisinage d'une composante fronti\`ere de cette classe
et expliquons enfin comment \'ecrire une sorte de
d\'eveloppement de Fourier et comment exprimer les coefficients de ce d\'eveloppement -- qui sont
des \'el\'ements du $H^1$ des courbes elliptiques fronti\`eres -- en terme du d\'eveloppement de
Fourier-Shimura (cf. \cite{Sh2} , \cite{Sh3}) de la forme de Picard de d\'epart. On peut donc lire sur ce d\'eveloppement
les propri\'et\'es de rationalit\'e de nos classes.  
\medskip
Il serait extr\^emement int\'eressant de pouvoir faire de m\^eme pour des classes de cohomologie 
de degr\'e 2 :  en effet, de m\^eme  que dans
\cite{Ca1}, on devrait obtenir  par cup-produit de formes de type holomorphe et anti-holomorphe toutes les
formes du type de Maass.  Si l'on pouvait d\'efinir des sortes de coefficients de Fourier pour de telles
classes, on en d\'eduirait des renseignements tout-\`a-fait nouveaux sur leurs propri\'et\'es arithm\'etiques.
Il s'agit toutefois d'un probl\`eme nettement plus difficile que pour le  $H^1$. On verra en effet
dans le cours de cet article \`a quoi ressemble un
voisinage de la courbe elliptique fronti\`ere : un  tel voisinage 
 n'est pas localement de Stein ;  le  $H^2$ d'un tel voisinage de notre courbe elliptique
fait intervenir, outre la g\'eom\'etrie de cette courbe, les groupes de $1$ - cohomologie des
voisinages des points de la fronti\`ere. Il s'agit donc d'un objet de nature nettement moins alg\`ebrique que le  $H^1$.
\medskip
Le plan de l'article est le suivant : au paragraphe 1, nous d\'efinissons les vari\'et\'es qui nous 
int\'eressent, et nous expliquons comment adapter les constructions de \cite{Ca2} afin de pouvoir
associer aux formes modulaires de Picard des classes de cohomologie automorphe.  Dans le second
paragraphe nous explicitons, dans le cas particulier consid\'er\'e ici, la construction de Kato et Usui, qui 
donne naissance \`a une compactification partielle de notre vari\'et\'e de Griffiths--Schmid.  Le troisi\`eme 
paragraphe est  consacr\'e \`a la d\'efinition des d\'eveloppements de Fourier-Jacobi tant dans le cas des 
formes modulaires de Picard que dans celui des classes de cohomologie automorphe sur les vari\'et\'es
ici consid\'er\'ees : le premier cas est bien classique, et on obtient des \'el\'ements du $H^0$ de courbes 
elliptiques CM, autrement dit des fonctions th\'eta. Dans le second cas nous obtenons comme "coefficients"
des classes dans le $H^1$ des m\^emes courbes - o\`u plut\^ot de leurs complexes conjugu\'ees.

Au  quatri\`eme paragraphe nous d\'efinissons une transformation cohomologique qui associe aux sections 
de faisceaux sur d'une courbe elliptique des classes de $1$--cohomologie sur la courbe conjugu\'ee et nous
\'etablissons des propri\'et\'es de rationalit\'e pour cette transformation. Finalement, nous v\'erifions au paragraphe 5
que la transformation cohomologique du \S1 est donn\'ee au niveau de chaque coefficient par celle du \S4,
ce qui prouve finalement les propri\'et\'es de rationalit\'e attendues.

\medskip
J'ai expos\'e une version pr\'eliminaire de ce travail lors d'un workshop "Arithmetic Groups and Automorphic Forms"
organis\'e au  d\'ebut 2002 par Joachim Schwermer \`a l'Erwin Schr\"odinger International Institute for Mathematical Physics de Vienne.  
Je remercie ici cet organisme ainsi que l'organisateur de la rencontre.

\vskip1cm

\head {\bf 1. Formes automorphes et cohomologie automorphe} \endhead
\bigskip
Soit $F \subset \Bbb C$ un corps quadratique imaginaire, muni d'un plongement 
complexe. On note $G =  SU
(2,1)$ le groupe  sp\'ecial unitaire (quasi-d\'eploy\'e) associ\'e \`a la forme hermitienne $H$
sur
$F^3$  donn\'ee par $H(x,y,t) = -x \bar t  + y \bar y - \bar x t$, o\`u $x\rightarrow
\bar x$ d\'esigne la conjugaison complexe sur $F$. Le groupe $G(\Bbb R)$ op\`ere naturellement
sur le plan projectif complexe $\Bbb P^2 (\Bbb C)$ avec deux orbites ouvertes 
qui sont respectivement la "boule unit\'e" ouverte $\Delta$ , constitu\'ee des points $p$
associ\'es aux vecteurs $v$ qui v\'erifient $H(v) < 0$, et le compl\'ementaire de 
la boule ferm\'ee $ \Delta^c$. On a repris
ici des notations de
\cite{Ca2} sauf en ce qui concerne la forme hermitienne qui \'etait diagonale dans (loc.
cit) et qu'il est plus pratique de consid\'erer maintenant comme nous le faisons ici, dans
le but d'\'etudier les ph\'enom\`enes qui se passent au voisinage des pointes. Une autre
diff\'erence essentielle avec \cite{Ca2} est que nous avions consid\'er\'e alors une forme
anisotrope de $SU(2,1)$. 
\medskip
Le domaine de p\'eriodes $\gW$ est le sous-espace
ouvert de l'ensemble des drapeaux $(p,L)$ constitu\'e de ceux tels que $p$
n'appartienne pas \`a la boule ferm\'ee $\Delta^c$ et que $L$ rencontre $\Delta$.
Si on note $\check \Omega$ l'ensemble de tous les drapeaux de  $\Bbb P^2 (\Bbb C)$, on voit qu'il 
y a 3 orbites ouvertes pour l'action de $G (\Bbb R)$ sur $\check \Omega$: d'une part $\gW$,
d'autre part deux orbites not\'ees $\bold X$ et $\bold Y$, 
en dualit\'e par rapport \`a la forme $H$: 
on a d\'esign\'e par
$\bold X$ (resp. $\bold Y$) 
le sous espace de l'espace des drapeaux $\check \bold \Omega$ constitu\'e des 
$(p,L)$ tels que $p \in \gD$ (resp. tels que $L \cap   \gD^c = \emptyset$).
\medskip
Les {\it vari\'et\'es de Griffiths-Schmid} (connexes) associ\'ees \`a cette situation sont les
quotients $\Gamma \backslash    \Omega$ pour $\Gamma$ un sous-groupe de congruence assez petit de
$G$. Ce sont des vari\'et\'es analytiques complexes non alg\'ebriques.
Contrairement \`a la situation \'etudi\'ee dans \cite{Ca2}, elles ne sont plus maintenant
compactes.
\medskip
Plus classiquement, on consid\`ere les {\it surfaces de Picard} $\Gamma \backslash    \Delta$.
L'espace 
$\bold X$ est fibr\'e en droites projectives au dessus de $\gD$, et il en est donc de m\^eme
pour les quotients $\gC \setminus \bold X$, au-dessus des surfaces  $\gC \setminus 
\gD$.  Quand \`a $\bold Y$ et ses quotients, ils s'identifient naturellement, sous la 
dualit\'e donn\'ee par la forme hermitienne $H$, aux conjugu\'es complexes de $\bold X$ et
de ses quotients.  Ils sont fibr\'es en droites projectives sur la surface de Picard 
 $\Gamma \backslash   \overline  \Delta$, o\`u $\overline \Delta$ d\'esigne l'ensemble,
dual de $\Delta$, constitu\'e 
des droites $L$ ext\'erieures \`a $\Delta^c$.
\medskip
D\'esignons d'autre part par $\Cal F_{a,b}$ la restriction  du faisceau 
$\Cal O (a) \otimes \Cal O (b) $ \`a $\check \Omega
\subset {\Bbb P}^2 (\Bbb C) \times {\Bbb P}^{ 2 \vee} (\Bbb C)$, ainsi qu'aux diff\'erents espaces
$\Omega$, $\bold X$, $\bold Y$; puisqu'il s'agit de faisceaux
\'equivariants , cela d\'efinit \'egalement des faisceaux,
not\'es de fa\c con identique, sur les  quotients $\Gamma \backslash    \gW$ ,
$\Gamma \backslash    \bold X$ , $\Gamma \backslash  \bold Y$. Nous appellerons {\it
formes modulaires de Picard} les sections, quand  elles existent, des faisceaux $\Cal F_{a,b}$ sur les
vari\'et\'es
$\gC \setminus
\bold X$ (resp. $\gC \setminus \bold Y$), munies de condditions de croissance convenable aux pointes.
 Pour le faisceau $\Cal F_{-k , 0}$ sur 
 $\gC \setminus \bold X$ (resp. $\Cal F_{0,-k}$ sur  $\gC \setminus \bold Y$), qui provient
de $\Gamma \backslash    \Delta$ (resp. $\Gamma \backslash   \overline \Delta$), on obtient
 ainsi la notion la plus usuelle de 'forme de Picard de poids $k$' ;  pour les autres on retrouve 
la d\'efinition plus habituelle en consid\'erant le fibr\'e vectoriel sur $\Gamma \backslash    \Delta$
(resp. $\Gamma \backslash  \overline   \Delta$) d\'efini comme l'image directe du faisceau constant
par la fibration en $\Bbb P^1$ dont il a \'et\'e question ci-dessus.

\bigskip

 Dans \cite{Ca2} nous avions d\'efini des transformations lin\'eaires:

$${\Cal P}: \ H^0(  {\bold X} \ , \ {\Cal F}_{a,b}) \longrightarrow
 H^1( { \Omega} \ , \ {\Cal F}_{-a-2, a+b+1}) \ , $$

$${\text et}\ \ \ {\Cal P}': \ H^0( {\bold Y} \ , \ {\Cal F}_{a,b}) \longrightarrow
 H^1( { \Omega} \ , \ {\Cal F}_{a+b+1,-b-2}) \ ,$$

\noindent
la premi\`ere \'etant injective pour $b \geq 0$ et $a+b \leq -2$, et la seconde
pour $a \geq 0$ et $a+b \leq -2$. Nous nous placerons toujours sous ces hypoth\`eses.
Nous appliquerons ces transformations aux sections $\gC$--invariantes qui correspondent
\`a des formes modulaires de Picard. On obtient alors comme image une {\it classe de 
cohomologie  automorphe}, c'est-\`a-dire invariante par $\gC$. Comme dans l'article 
pr\'ecit\'e, une telle classe provient d'une (unique) classe de cohomologie sur le
quotient $\gC \setminus \gW$ : une fa\c con de  voir ceci est de consid\'erer
 la suite exacte des
termes de bas degr\'es associ\'ee \`a la suite spectrale de Cartan--Leray
relative au quotient $\gC \setminus \gW$ :

$$ H^1 \bigl( \Gamma  ,  H^0 ({ \Omega}  ,  {\Cal F}_{a',b'}) \bigr)
 \rightarrow H^1({\Gamma}
\setminus { \Omega}  ,  {\Cal F}_{a',b'})
\rightarrow  H^1({ \Omega}  ,  {\Cal F}_{a',b'})^{\Gamma} \rightarrow
 H^2 \bigl( \Gamma  ,  H^0 ({ \Omega}  ,  {\Cal F}_{a',b'}) \bigr)  $$
avec $a'$ et $b'$ intervenant dans l'image de la transformation $\Cal P$ (resp. $\Cal P'$) ;
puis de remarquer que pour ces valeurs de $a'$ et $b'$, l'espace 
$H^0 ({ \Omega}  ,  {\Cal F}_{a',b'})$ est nul : en effet une telle section doit se
prolonger, d'apr\`es le lemme (4.2) de \cite{Ca2}, \`a l'espace de drapeaux $\check \gW$
tout entier, et le fait que $a+b+1 $ ( $=a'$ ou $b'$) soit n\'egatif entra\^\i ne qu'une
telle section ne peut qu'\^etre  nulle.

\medskip

Ce qui pr\'ec\`ede permet donc de d\'efinir des applications lin\'eaires injectives 
not\'ees encore $\Cal P$ et $\Cal P'$:

$${\Cal P}: \ H^0( \Gamma \setminus {\bold X} \ , \ {\Cal F}_{a,b}) \longrightarrow
 H^1( \Gamma \setminus{ \Omega} \ , \ {\Cal F}_{-a-2, a+b+1}) \ , $$

$${\Cal P}': \ H^0( \Gamma \setminus{\bold Y} \ , \ {\Cal F}_{a,b}) \longrightarrow
 H^1( \Gamma \setminus{ \Omega} \ , \ {\Cal F}_{a+b+1,-b-2}) \ .$$
\bigskip
Rappelons plus en d\'etail 
comment nous avions dans \cite{Ca2} d\'efini les transformations $\Cal P$ et $\Cal P'$
au niveau des espaces $\bold X$ , $\bold Y $, $\gW$ : nous avions utilis\'e la th\'eorie des
Eastwood-Gindikin-Wong (cf. \cite{EGW1} , \cite{EGW2}, \cite{Gi}), qui permet d'exprimer la cohomologie 
\`a partir d'une fibration \`a fibres contractiles et dont l'espace total est de 
Stein.  Dans notre cas, un tel espace est l'ensemble $\bold U$ constitu\'e des
 couples de drapeaux
 $(z,l;\xi, \alpha)$ v\'erifiant 
les conditions suivantes:

(i) les points $z$ et $\xi$ sont distincts et la droite $J$ qui les joint ne rencontre
pas $ \Delta^c$.

(ii) les droites $l$ et $\alpha$ sont distinctes et leur intersection $I$ appartient \`a
$\Delta$.
\medskip
 L'espace  $\bold U $ est de Stein et la projection sur le premier facteur $\pi : \bold U \rightarrow \gW$
est \`a fibres contractiles. Sous ces hypoth\`eses on a alors un isomorphisme entre la cohomologie 
de $ \Omega$ \`a valeurs dans un faisceau $\Cal F$ et la cohomologie du complexe 
$\Gamma ({\bold U} , \Omega_{\pi}^{\bullet} (\Cal F))$ des sections globales sur $\bold U$ du faisceau
 des diff\'erentielles relatives \`a valeurs dans $\Cal F$.  Nous avions alors d\'efini nos transformations 
dans ce cadre par des formules: 

$$ {\Cal P} (f) (z,l; \xi, \alpha) = f ( l\wedge \alpha , l) \ \alpha (z)^{-a} \ \omega_I $$
$$ {\Cal P}' (f') (z,l; \xi, \alpha) = f' ( z , z \wedge \xi) \ l(\xi)^{-b} \ \omega_J \  $$
avec $\omega_I \in \Gamma ({\bold U} , \Omega_{\pi}^1 ({\Cal F}_{-2,1}))$ et  
$ \omega _J \in 
\Gamma ({\bold U} , \Omega_{\pi}^1 ({\Cal F}_{1,-2})) $  des \'el\'ements canoniques  dont nous rappellerons 
plus bas (au \S5) la d\'efinition exacte.
\medskip
Ces m\^emes formules d\'efinissent \'egalement les transformations $\Cal P$ et $\Cal P'$ au
niveau des espaces quotients, une fois v\'erifi\'ee la proposition suivante (que nous avions 
prouv\'e dans (\cite{Ca2}  prop. 5.2) sous l'hypoth\`ese que $\Gamma $ \'etait co-compact) :

\proclaim{Proposition 1}
 L'espace $\Gamma \setminus {\bold U}$ est de Stein. 
\endproclaim
\demo{Preuve} La construction donn\'ee dans \cite{Ca2}  de fonctions qui
s\'eparent les points n'utilisait pas la co-compacit\'e de $\gC$ et reste
donc valide ici. 

 Par contre nous utilisions cette hypoth\`ese pour construire, \'etant
 donn\'ee une suite sans valeur d'adh\'erence $((I_n,l_n)  , (\xi_n, J_n))$
d'\'el\'ements de  $\Gamma \setminus {\bold U}$, 
une fonction $\Phi$ telle que $\Phi ((I_n,l_n)  , (\xi_n, J_n))$ ne
soit  pas born\'ee. Expliquons comment modifier l'argument afin qu'il s'applique au cas 
pr\'esent.
\medskip
 Tout d'abord on v\'erifie aussit\^ot que l'argument donn\'e dans (loc. cit.)
s'applique tel quel si  l'image de la  suite $I_n$ dans le quotient  $\gC \setminus
\Delta$, de m\^eme que celle  de la  suite $J_n$ dans le quotient  $\gC \setminus
\overline \Delta $, sont contenues dans des parties compactes. 
\medskip
  On peut donc supposer, quitte \`a permuter \'eventuellement le r\^ole de $\bold X$
et $\bold Y$, que $I_n$ admet une pointe comme valeur d'adh\'erence. Quitte \`a 
remplacer ensuite $\gC$ par un conjugu\'e, \`a extraire une sous-suite et \`a choisir
des repr\'esentants convenables, on se ram\`ene \`a supposer que $I_n$ est repr\'esent\'e par
 $\tilde I_n = \pmatrix
x_n \\ y_n \\ 1 \endpmatrix $, avec $y_n$ et $\Im (x_n)$ born\'es tandis que $\Re (x_n)$
tend vers $+\infty$.  On peut aussi supposer qu'un repr\'esentant  $\tilde J_n = (u_n, v_n , w_n)$
a \'et\'e choisi de telle sorte qu'il converge vers 
$\tilde J = (u, v , w)$, dont l'image $J$ peut appartenir \`a $\overline \gD$ ou \`a sa fronti\`ere suivant que
$\vert v \vert^2 - 2 \Re (u \bar v)$  est $< 0$ ou bien nul ; noter que dans le premier cas
$u$ est n\'ec\'essairement non nul.
\medskip
   Le cas le plus simple est le cas o\`u $J \in  \overline \gD$. Dans ce cas $f(\tilde I_n)$
converge vers la valeur de notre forme en la pointe consid\'er\'ee, $g(\tilde J_n)$ vers $g(\tilde J)$,
et $ \vert \tilde J_n (\tilde I_n) \vert  = \vert  u_n x_n + v_n y_n + w_n \vert $
tend vers $+\infty$. On voit alors que pour un choix convenable de $f$ et $g$, la suite
$$ \Phi_{f,g} ((I_n,l_n)  , (\xi_n, J_n))  = f (\tilde I_n)  g (\tilde J_n)  \tilde J_n (\tilde I_n)^{-a}$$
(notations de loc. cit avec $b = 0$ , $a \leq -2$) n'est pas born\'ee.
\medskip
  Le second cas est celui o\`u la suite $J_n$ converge vers un point de la fronti\`ere, mais o\`u
son image dans le quotient $\gC \setminus
\overline \Delta $ reste dans un compact. On peut alors trouver des scalaires $\gk_n$ ainsi que des
\'el\'ements $\gc_n \in \gC$ tels que $\gk_n  \tilde J_n \gc_n^{-1} = \tilde J'_n = \pmatrix u'_n & v'_n & w'_n
\endpmatrix$ converge vers  $\tilde J'=\pmatrix u' & v' & w' \endpmatrix$ repr\'esentant un 
\'el\'ement $J'$ de $\gD$. Dans ce cas 
$$\vert v'_n\vert^2 - 2 \Re (u'_n \bar w'_n) = \vert \gk_n \vert^2 (\vert v_n\vert^2 - 2 \Re (u_n \bar w_n))$$
tend vers un r\'eel $<0$ tandis que $\vert v_n\vert^2 - 2 \Re (u_n \bar w_n)$ tend vers $0$, de sorte que
$\vert \gk_n \vert$ tend vers l'infini.
\medskip
On a : $g (\tilde J'_n) = \gk_n^a g (\tilde J_n)$, de sorte que 
$$ \Phi_{f,g} ((I_n,l_n)  , (\xi_n, J_n))  = f ( \tilde I_n) )  g (\tilde J'_n)  
  (\gk_n \tilde J_n ( \tilde I_n))^{-a}
=f(\tilde I_n) g (\tilde J'_n) (\tilde J'_n \gc_n \tilde I_n)^{-a}.$$

Les vecteurs $\tilde I'_n = \gc_n \tilde I_n  = \pmatrix x'_n \\ y'_n \\ z'_n \endpmatrix $ repr\'esentent des points
de la boule unit\'e dont toutes les valeurs d'adh\'erence sont sur sa fronti\`ere (sans 
quoi leur projection sur le quotient  ne tendrait pas vers la fronti\`ere). Extrayant encore une sous-suite,
on peut supposer que pour certains scalaires $\gl_n$ le vecteur 
$\gl_n^{-1} \tilde I'_n$ converge vers  $\pmatrix x' \\ y' \\ z' \endpmatrix$
repr\'esentant un point fronti\`ere.  D'autre
part $\vert y'_n \vert ^2 - 2 \Re (x_n'\bar z'_n) = \vert y_n \vert ^2 - 2 \Re (x_n)$ tend vers $-\infty$
et donc $\vert \gl_n \vert$ tend vers $+ \infty$.
Enfin  $\gl_n^{-1} J'_n (I'_n)  $ tend vers  $u'x' + v'y' + w'z'$ qui est non nul, d'o\`u il d\'ecoule 
que la suite $ J'_n (I'_n)  $ n'est pas
born\'ee. Nous en d\'eduisons  comme pr\'ec\'edemment que, pour un choix convenable des fonctions $f$ et $g$, la
suite des 
$ \Phi_{f,g} ((I_n,l_n)  , (\xi_n, J_n))$ n'est pas born\'ee.
\medskip
 Il reste le cas o\`u l'image de $J_n$ converge aussi vers une pointe. Il existe alors $\gc \in G(\Bbb Q)$
tel que $J_n$ soit repr\'esent\'e par $\tilde J_n = \pmatrix 1 & v_n & w_n \endpmatrix \gc$ 
avec $v_n$ et $\Im w_n$
born\'es et
$\Re w_n$ tendant vers $+ \infty$. On a
$$f(\tilde I_n) g(\tilde J_n) \tilde J_n (\tilde I_n)^{-a} = 
f(\tilde I_n) g(\tilde J_n) (\pmatrix 1 & v_n & w_n \endpmatrix \gc
\tilde I_n)^{-a}.$$ 
Les $\gc \tilde I_n = = \pmatrix x'_n \\ y'_n \\ z'_n \endpmatrix$ repr\'esentent des points qui
convergent vers la fronti\`ere et comme plus haut pour les $\tilde I'_n$ on peut trouver  $\gl_n$ tendant vers
l'infini et tel que
$\gl_n^{-1}
\gc \tilde I_n$ converge vers 
$\pmatrix x' \\ y' \\ z' \endpmatrix$ repr\'esentant un point fronti\`ere. L'expression
$$\gl_n^{-1} \pmatrix 1 & v_n & w_n \endpmatrix \gc I_n
= \gl_n^{-1} ( x'_n + v_n y'_n + w_n z'n) $$ 
tend vers l'infini (si $z' \neq 0$ ), ou vers  $x' \neq 0$ dans le cas o\`u $z' = y' = 0$. Dans un 
cas comme dans l'autre, on voit que $\pmatrix 1 & v_n & w_n \endpmatrix \gc I_n$ n'est pas born\'ee
et on conclut comme dans les cas pr\'ec\'edents.
\enddemo
\bigskip
\noindent
{\bf Remarque} : Dans \cite{Ca2} nous avions donn\'e une expression des 
transformations $\Cal P$ et $\Cal P'$ en termes de cohomologie de Dolbeault et
qui ne fait pas appel \`a la th\'eorie de Gidinkin. Une telle expression est encore valide ici.
\medskip \noindent
{\bf Remarque} :  Dans (loc.cit.) nous avions montr\'e que les transformations 
 $\Cal P$ et $\Cal P'$  sont aussi surjectives.  Ceci est probablement encore vrai ici
\`a condition de se limiter aux sous-espaces constitu\'e des formes (resp. des classes)
paraboliques. Tout le probl\`eme est de d\'efinir une notion convenable de parabolicit\'e
pour les classes de cohomologie automorphe. On peut le faire sans difficult\'e d'un point de 
vue analytique, en exprimant de fa\c con habituelle  la cohomologie automorphe comme 
$(\Cal P, K)$--cohomologie d'un espace de formes automorphes puis en se limitant \`a la 
partie parabolique de ce dernier. Avec une telle d\'efinition, la surjectivit\'e des trandsformations
 $\Cal P$ et $\Cal P'$  est essentiellement triviale. Il serait int\'eressant d'avoir une notion  plus g\'eom\'etrique
analogue \`a ce que fait Harris (\cite{Ha}) dans le cas des vari\'et\'es de Shimura.   Malheureusement les compactifications
que nous allons utiliser dans cette article semblent trop partielles pour pouvoir produire une telle d\'efinition.

\vskip1cm

\head {\bf 2. Compactification de Kato -- Usui} \endhead
\bigskip
L'article \cite{K-U} construit des 'compactifications partielles' -- \`a vrai dire en
g\'en\'eral bien loin d'\^etre compactes -- des espaces classifiants de structures de Hodge
polaris\'ees. Cette m\^eme construction vaut encore si l'on ajoute des donn\'ees 
suppl\'e\-mentaires (actions). Nous allons appliquer cela \`a la vari\'et\'e de 
Griffiths--Schmid que nous consid\'erons. Comme la construction donn\'ee dans \cite{K-U}
  s'exprime en termes  de structures
de Hodge, nous allons commencer par expliquer bri\`evement (en suivant
\cite{De}) comment $\Omega$ param\`etre 
certaines structures de ce type. Cette description ne jouera dans la suite qu'un
r\^ole assez auxiliaire.
\medskip \noindent
{\bf (2.1)} Notons $W = F^3$,  muni de la forme hermitienne $H$, et $V$ 
le m\^eme espace apr\`es restriction
des scalaires \`a $\Bbb Q$. D\'esignons par $\Psi$ la forme altern\'ee sur $V$ obtenue comme
l'oppos\'ee de la partie imaginaire de $H$. Alors il revient au m\^eme de parler du groupe
unitaire (resp. des similitudes unitaires) de $(W,H)$, ou bien du 
groupe des \'el\'ements qui commutent \`a l'action de $F$ et
qui appartiennent au groupe symplectique (resp. des similitudes
symplectiques) de $(V,\Psi)$.

L'application $w \otimes \lambda \rightarrow (w \lambda, w \bar \lambda)$ identifie
$V \otimes F$ \`a la somme $W \oplus \overline W$ de $W$ et de son conjugu\'e $\overline W$
(le m\^eme espace mais muni de l'action conjugu\'ee de $F$). Avec cette identification
la forme $\Psi_F$ d\'eduite de $\Psi$ par extension des scalaires s'exprime comme il suit:
$$\Psi_F (v_1 \oplus \bar v_2 \ , \   w_1 \oplus \bar w_2) = {i \over 2} \big(
H(v_1 \ , \ \bar w_2) - H(w_1 \ , \ \bar v_2) \big) .$$

De m\^eme, on a $V_\Bbb R$, isomorphe \`a $W_\Bbb C = \Bbb C^3$ dont le complexifi\'e
s'identifie 
\`a  $W_\Bbb C \oplus \overline W_\Bbb C$. Choisissons une base $(e_1, e_2, e_3)$ de 
$W_\Bbb C$,  orthogonale et telle que $H(e_1) = H(e_2) = 1 = - H(e_3)$. Notons 
$\bar e_i$ les m\^emes \'el\'ements $e_i$, mais vus dans $\overline W_\Bbb C$, de sorte que 
l'on a : 
$$\Psi_F (e_i, \bar e_j ) = 0 \  {\text si }\  i \neq j$$ 
$$\Psi_F (e_1, \bar e_1 ) = \Psi_F (e_2, \bar e_2 ) =1 \ \Psi_F (e_3, \bar e_3 ) = -1 $$

Pour $z
\in
\Bbb C^*$ consid\'erons alors la similitude unitaire $h(z)$ de $W_\Bbb C$
dont la matrice dans cette base
s'\'ecrit $ \operatorname{diag} 
(z^{-p_1} {\bar z}^{-q_1}, z^{-p_2} {\bar z}^{-q_2},z^{-p_3} {\bar z}^{-q_3})$, avec six entiers
qui v\'erifient:

$p_1 + q_1 = p_2 + q_2 = p_3 + q_3 = w = -1$

$ - p_1 + q_1 \equiv - p_2 + q_2 \equiv 1 \mod 4$
et $ - p_3 + q_3 \equiv - 1 \mod 4$,

\medskip \noindent
cette seconde condition signifiant que $h(i) =  \operatorname{diag} (i,i,-i)$ d\'efinit 
une involution de Cartan ; elle entra\^\i ne d'autre part, compte tenu de la premi\`ere,
le fait que $p_1$ et $p_2$ sont pairs et $p_3$ impair (d'o\`u $q_1$ et $q_2$ impairs,
$q_3$ pair). De plus,
on a alors le fait que
$
\Psi (x, h(i) y)$ est une forme sym\'etrique d\'efinie positive sur $V_\Bbb R = W_\Bbb C$.
Finalement, prenant le
$\operatorname{F}^0$ de la structure de Hodge d\'efinie par $h$ sur l'alg\`ebre de Lie
complexifi\'ee de $G$, on trouve une sous-alg\`ebre parabolique qui co\"\i ncide avec
l'alg\`ebre de Lie du sous-groupe de Borel not\'e $B$ dans \cite{Ca1} { \bf (3.1)} si et
seulement si la condition suivante est satisfaite: $ p_1 > p_3 > p_2$. 
\medskip
Il y a de nombreux choix possibles de tels entiers, donnant lieu \`a des types diff\'erents de 
structures de Hodge (suivant la position relative des $p_i$ et 
des $q_j$) mais aboutissant finalement \`a la m\^eme
compactification. Prenons pour fixer  les id\'ees 
$p_2 = q_3 = -1$ ; $q_2 = p_3 = 0$ ; $p_1 = 1$ ; $q_1 = -2$.
 Dans ce cas la structure 
de Hodge associ\'ee \`a $h$ est donn\'ee par :
\medskip
$V^{1 , -2} = \langle e_1 \rangle \ ,\ \ \
V^{-1 , 0}  = \langle e_2 , \bar e_3 \rangle  \ , \ \ \
V^{-2 , 1}  = \langle \bar e_1 \rangle  \ , \ \ \
V^{0 , -1}  = \langle e_3 , \bar e_2 \rangle$ .

\medskip
D'autre part le Drapeau $(p,L)$ correspondant \`a $h$, c'est-\`a-dire celui fix\'e par $B$, est donn\'e
par: $p = \langle e_1 \rangle$ et $L = \langle e_1 , e_3 \rangle$. 
La filtration de Hodge sur $V$ associ\'ee est donn\'ee par :
\medskip
  $\operatorname{F}^0 = \{ 0 \}$ ;  $\operatorname{F}^1 = \langle e_1 \rangle $ ;   $\operatorname{F}^0 =
\langle e_1 , e_3 , \bar e_2
\rangle $ ;  $\operatorname{F}^{-1} = \langle e_1 , e_2 , e_3 , \bar e_2 , \bar e_3\rangle $ ;
$\operatorname{F}^{-2} = V$.
\medskip
Autrement dit on obtient comme filtration :

$$\operatorname{F}^1 = p \subset  \  \operatorname{F}^0 = L \oplus L^{\perp} \subset \operatorname{F}^{-1} =
W_\Bbb C \oplus p^{\perp} $$
(nous avons  d\'enot\'e ici -- un peu abusivement -- par les m\^emes 
notations le point $p$ (resp. la droite $L$) et les sous-espaces vectoriels de $W_\Bbb C$
correspondants).
\medskip
Prenant les conjugu\'es de $h$, lesquels correspondent aux diff\'erents points de $\gW$,  on obtient 
 des structures de Hodge  du type pr\'ec\'edent.
Ainsi $\gW$ classifie de telles structures (polaris\'ees, munies d'une action de $F$).
A chaque $(p,L) \in \gW $ correspond  une structure de Hodge, telle que la filtration de Hodge associ\'ee reste
donn\'ee par les m\^emes expressions que ci-dessus.
\bigskip \noindent
{\bf (2.2)} {\bf Donn\'ees combinatoires. }
\medskip
L'ingr\'edient  de base de la construction de \cite {K-U} consiste en la donn\'ee
d'un \'eventail $\Sigma$ constitu\'e de c\^ones $\sigma$ dans l'alg\`ebre de Lie 
sur $\Bbb Q$ de
$G$. Chacun de ces c\^ones doit \^etre engendr\'e par
des \'el\'ements nilpotents commutant entre eux. A une conjugaison pr\`es on voit que l'on
peut se ramener \`a supposer qu'un tel c\^one est contenu dans la sous-alg\`ebre, 
constitu\'ee des 
matrices triangulaires sup\'erieures :

$$\Cal N =  \Biggl\{ \ \pmatrix 0 & \ga  & \gb \\ 0 & 0 & \bar \ga \\ 0 & 0 & 0 \endpmatrix \ \  
/ \ \ \
  \ \ \gb + \bar \gb = 0\  \Biggr\} $$ 
et associ\'ee au sous-groupe:
$$ V =  \Biggl\{ \ \pmatrix 1 & \ga & \gb \\ 0 & 1 & \bar \ga \\ 0 & 0 & 1 \endpmatrix  \ \ /  \ \    
\ \   \ \ \gb + \bar \gb = \gaÊ\bar \ga  \  \Biggr\}
$$

Pour chaque tel c\^one $\gs$ on ajoute comme composante fronti\`ere
l'ensemble cons\-titu\'e des orbites $\gs$ - nilpotentes : ce sont les $\exp (\gs_\Bbb
C)$-orbites  $\exp (\gs_\Bbb C)(X)$
dans la vari\'et\'e de drapeaux $\check \gW$ qui  v\'erifient une condition de
transversalit\'e de
Griffiths et une condition de positivit\'e. La condition de positivit\'e
 exprime le fait que, pour $N_j$ des \'el\'ements du c\^one $\gs$ et pour des
r\'eels $Y_j$ tous assez grands, $ \exp (\sum i Y_j N_j )(X)$ 
appartient \`a $\Omega$. Celle de transversalit\'e est que, pour $N \in \gs$
et la filtration 
du type  Hodge qui continue \`a \^etre associ\'ee aux points de $\check \gW$, on ait 
$N \operatorname{F}^p \subset  \operatorname{F}^{p-1}$. La relation explicit\'ee ci-dessus entre drapeaux
$(p,L)$ et filtration de Hodge montre que cette condition est \'equivalente \`a :
 $N p \subset L$. 
\medskip
 Remarquons que, si $\gs$ est un c\^one de dimension 2, engendr\'e par deux 
\'el\'ements (non nuls, non proportionnels, commutant entre eux) :

$$N_1 = \pmatrix 0 & \ga_1  & \gb_1 \\ 0 & 0 &  \bar \ga_1 \\ 0 & 0 & 0
\endpmatrix \ \  {\text et}Ê\ \ 
 N_2 = \pmatrix 0 & \ga_2  & \gb_2 \\ 0 & 0 &  \bar \ga_2 \\
0 & 0 & 0
\endpmatrix$$
alors il n'y a aucune  orbite $\gs$-nilpotente. En effet si $X = (p,L)$ appartenait
\`a une telle orbite, avec $p =  \pmatrix x \\ y \\ t \endpmatrix$, on aurait $t \neq0$
par la condition de positivit\'e ;  la condition de 
transversalit\'e entra\^\i nerait que les points d\'efinis par :

$$N_1 p =  \pmatrix \ga_1 y + \gb_1 t \\ \bar \ga_1 t \\ 0 \endpmatrix  \ \  
{\text et} \ \ 
\ \ \ N_2 p =  \pmatrix \ga_2 y + \gb_2 t \\ \bar \ga_2 t \\ 0 \endpmatrix$$
appartiennent \`a la droite $L$. Comme ils appartiennent tous deux \`a  la "droite \`a
l'infini" d'\'equation $t=0$, ils sont donc confondus, ce qui s'\'ecrit : 
$$(\ga_1 \bar \ga_2 - \bar \ga_1 \ga_2) ty + (\bar \ga_2 \gb_1 - \bar \ga_1 \gb_2) t^2 = 0 ;$$
d'autre part la condition de commutation entre $N_1$ et $N_2$ est que 
$\ga_1 \bar \ga_2 - \bar \ga_1 \ga_2 = 0$; il en r\'esulte que $\bar \ga_2 \gb_1 - \bar \ga_1 \gb_2 = 0$,
ce qui contredit  la non-proportionalit\'e de $N_1$ et $N_2$. 
\medskip
{\it On peut donc se limiter \`a ne consid\'erer que des c\^ones de dimension $1$.}
Un tel c\^one est 
  engendr\'e par un \'el\'ement nilpotent $N$. On normalisera dans la suite ce nilpotent  comme 
dans \cite{K-U},
de telle sorte que $\exp (N)$ soit un g\'en\'erateur de $\gC(\gs) = \exp (\gs) \cap \gC$ . De tels $N$ sont de
deux types possibles suivant leur ordre de nilpotence ($2$ ou $3$). Les \'el\'ements nilpotents d'ordre
$3$ et qui appartiennent \`a $\Cal N$ sont ceux de la forme 
$ \pmatrix 0 & \ga  & \gb \\ 0 & 0 & \bar\ga \\ 0 & 0 & 0 \endpmatrix$ avec $\ga \neq 0$ ;
ceux d'ordre $2$ 
 s'\'ecrivent $ \pmatrix 0 & 0  & \gb \\ 0 & 0 & 0 \\ 0 & 0 & 0
\endpmatrix$. 

\bigskip  \noindent {\bf (2.3)} Nous allons d\'eterminer
pour chacun de ces types l'espace des orbites $N$--nilpotentes, c'est-\`a-dire la composante
fronti\`ere qui lui correspond.  Nous ne traitons le cas des nilpotents d'ordre 3 qu'\`a titre indicatif,
car seules joueront un r\^ole dans la suite de l'article les
composantes attach\'ees aux nilpotents d'ordre 2.

\smallpagebreak
{\bf (a) Cas d'un nilpotent d'ordre $3$ :} 
$N =  \pmatrix 0 & \ga  & \gb \\ 0 & 0 & \bar \ga \\ 0 & 0 & 0 \endpmatrix$.  Faisant agir \break
$\exp (i Y N) = \pmatrix 1 & i Y \ga   &  - {Y^2 \over 2}\ga \bar \ga + i Y \gb  \\ 
0 & 1 &  i Y \bar \ga  \\ 0 & 0 & 1 \endpmatrix $ sur un drapeau  $X = (p,L)$ 
avec \break $p =  \pmatrix x \\ y \\ t \endpmatrix$ et
$L = \pmatrix u & v  & w  \endpmatrix$
on obtient le drapeau $X' = (p',L')$ avec \linebreak $p' =   \pmatrix x' \\ y' \\ t'
\endpmatrix = \pmatrix x + i Y \ga  y
  - {Y^2 \over 2} \ga \bar \ga  t + i Y \gb t \\ y + i Y \bar \ga  t\\ t \endpmatrix$ et 
$L' = \pmatrix u' & v'  & w'  \endpmatrix $ \linebreak $ =\pmatrix u & - i Y \ga  u + v   &  -
{Y^2Ê\over 2}  \ga \bar \ga   u - i Y \gb u
 - i Y \bar \ga  v + w\endpmatrix$. La condition de positivit\'e est que, pour $Y$ r\'eel assez 
grand, $y' \bar y' - 2 \Re (x' \bar t')$ ainsi que 
 $v' \bar v' - 2 \Re (w' \bar u')$ soient tous deux $>0$. D'autre part 
on doit avoir que $N p = \pmatrix \ga y + \gb t\\ \bar \ga t \\ 0 \endpmatrix \in L$. On en d\'eduit  
 que $t \neq 0$ : sans quoi en effet, puisque $L$ ne peut pas \^etre la droite d'\'equation
$t=0$, on voit que $p$ et $Np$ seraient li\'es, ce qui entra\^\i ne  $y=0$ et conduit
finalement \`a une contradiction avec la condition de positivit\'e.
\medskip
Si $t \neq 0$ la condition de positivit\'e pour $p'$ est satisfaite ($Y$ grand). Celle pour $L'$
l'est aussi : c'est clair si $v \neq 0$ ; sinon on remarque que $u \neq 0$ : en effet,  $L$, passant par 
$Np$,  ne passe pas par le point $\pmatrix 1 \\ 0 \\ 0\endpmatrix$.

 En r\'esum\'e on doit partir d'un
drapeau $(p,L)$  avec  $p =  \pmatrix x \\ y \\ 1\endpmatrix$ et $L$ la droite joignant $p$
\`a $Np$.  La $\exp (\gs_\Bbb C)$-orbite correspondante se compose des 
\'el\'ements $(p',L')$ avec $p'$ donn\'e par la m\^eme formule que ci-dessus -- mais $Y$
variant dans $\Bbb C $ tout entier -- et $L'$ la	 droite correspondante, joignant $p'$
\`a $Np'$. Chaque orbite de ce type admet un unique repr\'esentant v\'erifiant 
$p' =  \pmatrix x \\ 0 \\ 1\endpmatrix$, ce qui permet de param\'etrer cet ensemble 
d'orbites par $x \in \Bbb C$.  Version duale : on part de la droite $L = (1, v, w)$ et 
$p$ le point d'intersection de $L$ et $L N = (0, \ga ,  \gb + \bar \ga v)$
\bigskip
{\bf (b) Cas d'un nilpotent d'ordre $2$ } :
Si $\gb_0 $ est le plus petit \'el\'ement imaginaire, de partie imaginaire positive,
tel que $\pmatrix 1 & 0 & \beta_0 \\ 0 & 1 & 0 \\ 0 & 0 & 1 \endpmatrix \in \gC$, on peut 
supposer que 
$N =  \pmatrix 0 & 0 &  \pm \beta_0 \\ 0 & 0 & 0 \\ 0 & 0 & 0 \endpmatrix$.
 Alors ${ \exp} (i Y N) = \pmatrix 1 & 0   & \pm i Y \gb_0  
\\  0 & 1 &  0  \\ 0 & 0 & 1 \endpmatrix $ transforme  $X = (p,L)$ en $X' = (p',L')$,
donn\'e par (avec les m\^emes notations que ci-dessus) :
$$x' = x \pm i Y \gb_0 t \ , \  y' = y \ , \  t' = t \ ; \
 u' = u \ , \  v' = v \ , \  w' = \mp i Y \gb_0 u + w \ .$$

La conditions de positivit\'e correspondante signifie que pour $Y$ r\'eel assez 
grand: 
\medskip

 $ Ê \hskip .9cm  y \bar y - 2 \Re (x \bar t) \mp 2 i Y \beta_0 t \bar t > 0$

  et  : $ v \bar v - 2 \Re (w \bar u) \pm 2 i Y \beta_0 u \bar u > 0$
\bigskip
{\bf Cas (b +)}. $N =  \pmatrix 0 & 0 & \beta_0 \\ 0 & 0 & 0 \\ 0 & 0 & 0 \endpmatrix$. 
Alors la premi\`ere des deux conditions pr\'ec\'e\-dentes est  
satisfaite si $t \neq 0$ ou si $t =0$ et $y \neq 0$ ,
tandis que la seconde l'est si et seulement si $u = 0$ et  
$ v \neq 0 $. Les conditions sont donc \'equivalentes  \`a ce  que la droite $L$ passe par le
point 
$p_{\infty}  =  \pmatrix 1 \\ 0 \\ 0 \endpmatrix$ (ce qui entra\^\i ne la
condition $N p \  (= p_{\infty}) \in L$) et soit distincte de la droite  
$L_{\infty} = \pmatrix 0 & 0  & 1  \endpmatrix$, et que $p$ soit un point de cette 
droite $L$ distinct de $p_{\infty}$. La $\gs_\Bbb C$--orbite correspondante se d\'ecrit ainsi :
$L$ est fixe et $p$ varie sur $L \setminus \lbrace p_{\infty} \rbrace$. Remarquer que l'on a finalement que 
$t \neq 0$.
\bigskip
{\bf Cas (b -)} .  $N =  \pmatrix 0 & 0 & -\beta_0 \\ 0 & 0 & 0 \\ 0 & 0 & 0 \endpmatrix$. Maintenant 
la premi\`ere condition entra\^\i ne $t = 0$ et $y \neq 0$ : $p$ varie sur $L_{\infty}$
priv\'ee de $p_{\infty}$, et $L$ est une droite passant par $p$ et distincte, en 
vertu de la seconde condition, de  $L_{\infty}$. L'orbite correspondante est telle
que $p$ est fixe et $L$ varie en passant par $p$ et en restant distincte de $L_{\infty}$. 
\medskip
Remarquer que l'on a sur $\gW$ et ses quotients une involution anti-holomorphe pro\-venant de
la forme hermitienne.  Cette involution pr\'eserve l'ensemble des orbites du type (a) et
\'echange celles de type (b +) et (b -).

\medpagebreak \noindent
{\bf (2.4). Compactification de \cite{K-U}}. Pour chaque c\^one $\gs$, on ajoute comme composante fronti\`ere
l'ensemble des orbites $\gs$-nilpotentes, quotient\'e par le stabi\-lisateur dans $\gC$ de $\gs$. Il est plus 
difficile de voir comment recoller cette composante \`a la vari\'et\'e $\gC \setminus \gW$. Dans le cas
d'un c\^one de dimension 1 engendr\'e par un \'el\'ement $N$, on introduit dans ce but l'espace :

\medskip
\noindent
$$E_\gs = \Biggl\{(\gth,X) \in \Bbb C \times \check \Omega \ \Biggm\vert \  \Bigl\{ \alignedat 2
\text{si} \ \gth \neq 0 \ ,  & \exp((\log (\gth) / 2 \pi i) N) X \in
\Omega\ ; \\ \text{si} \  \gth= 0  \ ,  &  \exp(\gs_\Bbb C ) X \text{ est   
 une  orbite} \ \gs-\text{ nilpotente.}  \endalignedat  \Biggr\} \ $$
\medskip
On obtient ensuite $\Gamma_\gs^{\text grp} \setminus\Omega_\gs$, quotient de $\Omega$
auquel on a accol\'e  la composante fronti\`ere associ\'ee \`a $\gs$ par le groupe engendr\'e
par $\Gamma (\gs)$ , comme le quotient
de $E_\gs$ par l'action de $\Bbb C$ donn\'ee par :

$$\gl . (\gth , X) = (\exp (2 \pi i \gl) \gth , \exp ( - \gl N) X) .$$

Il ne reste plus alors qu'\`a quotienter l'espace ainsi obtenu par le stabilisateur dans $\gC$ de $\gs$,
et ce quotient d\'ecrit un voisinage de la composante fronti\`ere consid\'er\'ee.
\bigskip
Mettons en pratique cette construction dans chacun des cas consid\'er\'es.

\medskip
{\bf  ( Cas a ) : } On peut supposer pour fixer les id\'ees et simplifier les notations que
$N =  \pmatrix 0 & 1  & 0 \\ 0 & 0 & 1 \\ 0 & 0 & 0 \endpmatrix$ . Alors 
$E_\gs$ se compose des triplets $(\gth, p, L)$ \ 
avec $p =  \pmatrix x \\ y \\ t \endpmatrix$ et
$L = \pmatrix u & v  & w  \endpmatrix$\ qui v\'erifient les conditions suivantes, o\`u l'on a
pos\'e $q = \log (\gth) / 2 \pi$:
\medskip

- Si $\gth \neq 0$ : 

$$ \aligned \vert y - i q t \vert^2 - 2 \Re \  ( (x - i q y - {q^2 \over 2} t) \ \bar
t )
\ >0  \\ {\text et : }Ê \  \  \  \vert v + i q u \vert^2 - 2 \Re \  ( (w + i q v - {q^2 \over2}
u)
\ \bar u )  >Ê \ 0 \endaligned $$

- Si $\gth = 0 $ : $t \neq 0$ et $L$ est la droite joignant $p$ \`a $Np = \pmatrix y \\ t \\ 0
\endpmatrix$

L'action de $\gl \in \Bbb C$ est donn\'ee par: 

$$\gl . (\gth \    ,  \ ( p , L)) = (\exp (2 \pi i \gl) \gth \ , \ (p' , L'))$$
avec $p' = \pmatrix x - \gl y + {\gl^2 \over 2} t \\ y - \gl t \\ t \endpmatrix$ et 
$L' = \pmatrix u & \gl u + v  & {\gl^2 \over 2} u + \gl v + w  \endpmatrix$. Comme on cherche 
\`a d\'ecrire un voisinage de la composante fronti\`ere on peut se placer dans
l'ouvert $E'_\gs$ de $E_\gs$ constitu\'e des points tels que $t\neq 0$, lequel est
stable par l'action de $\Bbb C$. Dans cet ouvert, l'ensemble des points pour lesquels
$y=0$ constitue un ensemble de repr\'esentants transverse aux orbites (un "slice"), de
sorte qu'on peut d\'ecrire un voisinage de la composante fronti\`ere comme l'ensemble
des points $(\gth \ , \pmatrix x \\ 0 \\ 1 \endpmatrix , \pmatrix u & v  & - u x   \endpmatrix)$
qui v\'erifient les conditions ci-dessus. En particulier, pour $\gth = 0$ on doit avoir 
$v = 0$ et par suite $u \neq 0$. De sorte que l'on peut se restreindre un peu plus, et se
placer dans l'ouvert o\`u $uÊ\neq 0$. Notre nouveau voisinage est l'ensemble des
$(\gth \ , \pmatrix x \\ 0 \\ 1 \endpmatrix , \pmatrix 1 & v  & -  x   \endpmatrix)$, pour
lequel les conditions pr\'ec\'edentes se r\'ecrivent : 

\medskip
- Si $\gth \neq 0$ : 

$$2 (\Re (q))^2 - 2 \Re (x) > 0 \  \ {\text et} \ \ 2 (\Re (q))^2 + \vert v \vert^2 + 2 \Re (x) + 4 
\Re (q) \Im (v) > 0$$

- Si $\gth = 0 $ : $v = 0$.
\medskip
Au voisinage de $\gth = 0$ les in\'egalit\'es pr\'ec\'edentes sont automatiquement
satisfaites  ($\Re (q)$ tend vers
$-\infty$),  de sorte qu'un voisinage de cette composante fronti\`ere peut se d\'ecrire 
localement comme la droite ${\lbrace (0 , x , 0) \rbrace}$  ajout\'ee \`a l'ensemble des 
\lbrace $ (\gth , x , v ) \in \Bbb C^* \times \Bbb C \times \Bbb C \rbrace$. Il reste ensuite
\`a passer au quotient par le stabilisateur dans $\Gamma$ de notre c\^one $\gs$. On voit que ce
stabilisateur est constitu\'e (en supposant $\gC$ assez
petit) des matrices de la forme $\pmatrix 1 & \ga' & \gb' \\ 0 & 1 &  \ga' \\ 0 & 0 & 1
\endpmatrix$ avec $\ga' \in \Bbb Z$ et $\gb'$ entier de $F$ v\'erifiant $\gb' + \bar \gb' = (\ga')^2$ , ces
deux
\'el\'ements \'etant soumis \`a des conditions de congruence.  On v\'erifie qu'une telle matrice
envoie $( \gth , x , v ) $ sur  $( \gth , x + i \Im (\gb') , v ) $. En particulier la composante
fronti\`ere ajout\'ee s'identifie au quotient de $\Bbb C$ par un r\'eseau de $i \Bbb R$,
et donc \`a une copie de $\Bbb C^*$.

\bigskip
{\bf ( Cas b ) :} 
$N =  \pmatrix 0 & 0 & \pm \beta_0 \\ 0 & 0 & 0 \\ 0 & 0 & 0 \endpmatrix$ . Dans ce cas la
condition lorsque $\gth \neq 0$ \break s'\'ecrit :

$$ \aligned \vert y \vert^2 - 2 \Re \ (x \bar t)  \pm 2 (\Re (q)) i \gb_0 \vert t \vert^2  > 0
  \\\text{ et} \ Ê\ \vert v  \vert^2 - 2 \Re \ (w \bar u)  \mp 2 (\Re (q)) i \gb_0  \vert u
\vert^2  > 0 \endaligned $$
\noindent
L'action de $\gl \in \Bbb C$ envoie $(\gth , (p , L )) $ sur $((\exp (2 \pi i \gl) \gth ,
(p' , L'))$ avec \break $p' = \pmatrix x \mp \gl \gb_0 t  \\ y  \\ t \endpmatrix$ et
$L' = \pmatrix u &  v  & \pm \gl \gb_0 u  + w  \endpmatrix$.

\medskip
{\bf (Cas b +)} :  On a explicit\'e ci-dessus la condition que l'on a 
pour $\gth = 0$ :
en particulier, on a alors $t \neq 0$ et $v \neq 0$
. On peut se placer dans un voisinage $E'_\gs$
de la  composante fronti\`ere o\`u ces conditions sont encore remplies. Ce voisinage est stable
par l'action de $\Bbb C$. On obtient un syst\`eme de 
repr\'esentants transverse aux orbites dans $E'_\gs$ en faisant : $x = 0$ ; d'autre part, on peut
prendre $t=1$ et $v = 1$ , et donc $w = - y$. Dans $E'_\gs$ la premi\`ere moiti\'e des conditions
relatives au cas
$\gth \neq 0$ est
automatique d\`es que $\gth$ est assez petit. Il en r\'esulte  que l'on peut d\'ecrire
localement un voisinage de la composante fronti\`ere $\lbrace (\gth , y , u ) =  (0 , y , 0
)\rbrace$ en ajoutant \`a cette derni\`ere l'ensemble des $(\gth , y , u )$ qui v\'erifient
l'in\'egalit\'e :

$$ 1 + 2 \Re \ (\bar y u) - 2 (\Re (q)) i \gb_0  \vert u \vert^2  > 0$$

ou, si l'on pr\'ef\`ere: 

$$\log \vert \gth \vert >  \pi {1 + 2 \Re (\bar y u )\over i \gb_0 \vert u \vert^2}  \ \ . $$
On doit ensuite passer au quotient par $\gC \cap V$. Une matrice 
 $\pmatrix 1 & \ga & \gb \\ 0 & 1 &  \bar \ga \\ 0 & 0 & 1\endpmatrix$ ( avec $\gb + \bar \gb = \ga \bar \ga$)
envoie $(\gth , y, u)$ sur le point $(\gth , p' , L')$ o\`u
 $p' = \pmatrix \ga y + \gb \\ y + \bar \ga \\ 1 \endpmatrix$ et \break
$L' = \pmatrix u &  1 - \ga u  & u \bar \gb - \bar \ga -  y \endpmatrix$ et ce dernier 
est \'equivalent par l'action de $\gl = {1 \over \gb_0}(\ga y + \gb)$ \`a 
$\exp ({2 \pi i \over \gb_0} (\ga y + \gb)) \ \gth\ , \ p'' \ , \ L'')$ avec
$p'' = \pmatrix 0 \\ y + \bar \ga \\1 \endpmatrix$ et $$ L" = \pmatrix
u & 1 - \ga u & \ga yu + \gb u  + u \bar \gb - \bar \ga -  y \endpmatrix 
= (1 - \ga u ) \pmatrix {u \over 1 - \ga u} & 1 & -\bar \ga -  y \endpmatrix .$$
\medskip
 En d\'efinitive
$(\gth , y, u)$ est envoy\'e sur $ (\exp ({2 \pi i  \over \gb_0} (\ga y + \gb)) \ \gth \  ,  \ y + \bar \ga ,
{u \over 1 - \ga u} ).$  En parti\-culier, la composante fronti\`ere est $\Bbb C$ quotient\'ee
par un r\'eseau (un ordre $\Cal R$ de $F$ plong\'e par $\ga \rightarrow \bar \ga$), et donc une courbe
elliptique CM que nous noterons $\Cal E$.  Noter que le voisinage $E_\gs$  que nous avions choisi n'est pas
stable sous
$\gC \cap V$ (la condition $v \neq 0 $ n'est pas conserv\'ee) ce qui explique le d\'enominateur 
dans la formule pr\'ec\'edente (laquelle d\'ecrit n\'eanmoins la structure locale au
voisinage de notre composante fronti\`ere).

\medskip
{\bf (Cas b -) :} Dans ce cas, pour $\gth = 0$ on a  $y \neq 0$ et $u \neq 0$
 de sorte qu'on peut se placer dans le voisinage $E'_\gs$ d\'efini par ces conditions.
Alors un syst\`eme transverse de repr\'esentants est constitu\'e des points qui v\'erifient
$w = 0$. On peut prendre $y = u = 1$, de sorte que $v = - x$ . Cette fois c'est la seconde 
des deux in\'egalit\'es qui est automatique pour $\gth$ assez petit. On d\'ecrit localement 
notre compactification partielle
 en ajoutant \`a la droite$(\gth , x , t) = (0 , x , 0)$ l'ensemble des points
d\'efini par :

$$\log \vert \gth \vert >  \pi {1 - 2 \Re (\bar x t )\over i \gb_0 \vert t \vert^2} \ \  . $$

Maintenant l'action d'une matrice comme ci-dessus envoie notre point sur \linebreak
$(\gth , p' , L')$ avec
 $p' = \pmatrix x + \ga + \gb t  \\ 1 + \bar \ga t \\ t \endpmatrix$ et
$L' = \pmatrix 1 &  -x - \ga   &  \bar \gb + \bar \ga  x \endpmatrix$ ; faisant agir \break
$\gl =  {\bar \gb + \bar \ga x \over \gb_0}$ on obtient en d\'efinitive le point associ\'e 
\`a : \hfill \break \noindent$ (\exp ( {2 \pi i  \over \gb_0}(\bar \ga x + \bar \gb)) \ \gth \  , \ x + \ga  ,
{t \over 1 + \bar \ga t} ).$  La fronti\`ere est la courbe elliptique quotient par le m\^eme
r\'eseau que ci-dessus, mais plong\'e par $\ga \rightarrow \ga$, autrement dit la complexe 
conjugu\'ee $\Cal E'$ de la pr\'ec\'edente.

\medpagebreak \noindent
{\bf (2.5)  Remarques.} La structure des compactifications au voisinage des composantes de type (b)
est compliqu\'ee ; en particulier un voisinage \'epoint\'e d'un point-fronti\`ere n'est pas de Stein. Ceci en contraste 
avec les composantes du type (a) qui ont une structure beaucoup plus simple, comme expliqu\'e ci-dessus.
\smallskip
Pour voir cela, par exemple dans le cas (b +), on peut se placer dans un voisinage du point de la 
courbe fronti\`ere image de $y = 0$ et 
consid\'erer la section transverse d\'efinie par $y = 0$.  C'est donc l'intersection d'un  voisinage du point $(0,0)$ avec 
l'ensemble constitu\'e des $(\gth , u)$
qui v\'erifient 
$$\log \vert \gth \vert >  \pi {1 \over i \gb_0 \vert u \vert^2}  .$$

Posant $\Theta = \log \gth$ et $U = \log u$, on voit que l'in\'egalit\'e pr\'ec\'edente
se r\'ecrit : 
$$ \Re \Theta >  {\pi \over i \gb_0} \exp ( - 2 \Re U)$$
L'ensemble des tels couples $(\Theta , U)$ constitue donc un domaine tube au dessus d'un 
ouvert concave de $\Bbb R^2$ (rappelons en effet que $i \gb_0$ est un r\'eel n\'egatif).  Il n'est 
donc pas de Stein, non plus que l'ensemble des $(\gth, u)$ dont il constitue un rev\^etement.
De m\^eme pour les voisinages de $(0,0)$, dont en obtient les rev\^etements d'un  syst\`eme fondamental en tronquant
par les conditions $\Theta  < -A$  et  $U < -A$ pour $A$ un r\'eel positif assez grand. 
\medskip
Il serait tr\`es int\'eressant de comprendre la cohomologie d'un tel voisinage. Com\-me nous l'avons
d\'ej\`a dit dans l'introduction, cela permettrait d'associer des coefficients de Fourier-Jacobi aux classes de
cohomologie de degr\'e 2 et d'en d\'eduire des renseignements arithm\'etiques nouveaux sur les formes
modulaires du type de Maass.  Faute d'une telle compr\'ehension, nous allons dans la suite nous borner 
\`a d\'efinir de tels coefficients dans le cas des classes de degr\'e 1, ce qui n'utilise que la partie la plus 
alg\'ebrique de la cohomologie.

\vskip1cm
\head {\bf 3. D\'eveloppements \`a la fronti\`ere des formes modulaires et des classes de cohomologie automorphe} \endhead
\bigskip \noindent
{\bf (3.1)} {\bf Le cas classique : D\'eveloppement des formes modulaires de \break Picard. }
\medskip
Pour simplifier, nous nous bornerons ici \`a ne consid\'erer que les faisceaux $\Cal F_{-k , 0}$ sur 
 $\gC \setminus \bold X$ (resp. $\Cal F_{0,-k}$ sur  $\gC \setminus \bold Y$), qui proviennent
des surfaces de Picard $\Gamma \backslash    \Delta$ (resp. $\Gamma \backslash   \overline \Delta$).
Nous ne consid\'ererons de plus que des sections  qui s'annulent aux pointes. 
De telles
 sections correspondent, comme on l'a d\'ej\`a fait remarquer, \`a des  'formes de Picard paraboliques de poids $k$' .

Expliquons comment \'ecrire le d\'eveloppement de Fourier -- Shimura d'une telle forme ; les coefficients 
qui apparaissent sont des fonctions th\'eta, associ\'ees \`a certaines courbes elliptiques C M. 
G\'eom\'etriquement, cela correspond au fait que les surfaces de Picard admettent des compactifications 
lisses, dont la fronti\`ere est compos\'ee de telles courbes : voir \cite{Lar} pour la construction d'un mod\`ele 
arithm\'etique d'une telle compactification.

Nous suivrons ici le point de vue plus concret de  Shimura (\cite{Sh1} et \cite{Sh2}). On peut voir $\gD$ 
comme l'ouvert du plan 
constitu\'e des points
$(x,y)
\in
\Bbb C^2$ qui
 v\'erifient l'in\'egalit\'e $ 2 \Re (x) > \vert y \vert^2$.  Nous allons d\'efinir le
d\'eveloppement de Fourier de  formes modulaires, en une pointe associ\'ee au parabolique normalisateur
du sous-groupe $V$ d\'ej\`a consid\'er\'e plus haut. Une forme de Picard s'identifie \`a une
fonction $f(x,y)$ avec $(x,y)$ les coordonn\'ees d'un point de $\gD$, qui doit 
v\'erifier la loi de transformation habituelle, et en
particulier \^etre invariante sous le groupe $\gC \cap V$. Autrement dit, pour 
$   \pmatrix 1 & \ga & \gb \\ 0 & 1 & \bar \ga \\ 0 & 0 & 1 \endpmatrix  \in \gC \cap V $, on doit avoir :
$f(x + \ga y +\gb \ ,  \ y + \bar \ga) = f(x,y)$.

 Consid\'erant le cas o\`u $\ga=0$, et donc $\gb$ imaginaire 
multiple de $\gb_0$, on voit que  ce qui pr\'ec\`ede et la condition de cuspidalit\'e entra\^\i nent que $f$ admet un
d\'eveloppement:

$$f(x,y) = \sum_{r \in \Bbb N^*} g_r (y) \exp ( -{2 \pi  i r \over \gb_0} x) $$
avec  $g_r$  satisfaisant \`a l'\'equation
fonctionnelle :
$$g_r (y + \bar \ga) = g_r (y) \exp ({2 \pi  r i \over \gb_0} (\ga y + \gb)) .$$
Il s'agit donc de fonctions th\'eta, associ\'ees \`a la courbe elliptique $\Cal E$
d\'ej\`a rencontr\'ee plus haut. La relation ci-dessus d\'efinit un fibr\'e en
droites $\Cal L_r$ sur $\Cal E$ dont les $g_r$ constituent des sections. Plus pr\'ecis\'ement on peut
r\'ecrire  la relation pr\'ec\'edente sous la forme :
$$g_r (y + \bar \ga) = g_r (y) \chi (\bar \ga) \exp ({2 \pi  r i \over \gb_0} (\ga y + {1 \over 2} \ga \bar \ga)) $$
o\`u le multiplicateur $\chi$ est d\'efini par
$\chi (\bar \ga) = \exp ({2 \pi  r i \Im \gb \over \Im \gb_0})$ (en effet $\gb$ est d\'etermin\'e par 
$\ga$  \`a un multiple entier pr\`es de $\gb_0$) ; $\chi (\bar \ga)$ est une racine de l'unit\'e 
parce que $\Im \gb \over \Im \gb_0$ est un rationnel.
\medskip
De fa\c con duale, on peut consid\'erer des formes paraboliques dans $H^0 ( \Gamma \setminus{\bold Y} \ , \ {\Cal
F}_{0,b})$. Une telle forme $f'$ ne d\'epend que de $L$,  de coordonn\'ees $(1,v,w)$,  et satisfait \`a
la loi de transformation :
$f'(v-\ga , w- \bar \ga v + \bar \gb) = f' (v, w)$. 
Il en r\'esulte qu'elle admet un d\'eveloppement :
$$f'(v,w) = \sum_{r \in \Bbb N^*} g'_r (v) \exp ( -{2 \pi  i r \over \gb_0} w) $$
o\`u les $g'_r$ satisfont \`a
$$g'_r (v + \ga ) = g'_r (v) \exp (2 {\pi  r i \over \gb_0} (\bar \ga v + \gb)) $$
autrement dit ce sont des fonctions th\'eta relatives \`a la courbe elliptique conjugu\'ee
$\Cal E'$ et au faisceau conjugu\'e  $\Cal L'_r$.

\bigskip \noindent
{\bf (3.2)} {\bf D\'eveloppement de Fourier au voisinage de la fronti\`ere de Kato-Usui }

\medskip
Commen\c cons par remarquer que les faisceaux $\Cal F_{a,b}$ sont triviaux au voisinage 
des composantes fronti\`eres de type {\bf (b)}  : pour le type  {\bf (b +)} cela r\'esulte 
du fait, d\'ej\`a remarqu\'e en {\bf (2.4)}, que dans un tel voisinage on peut
supposer $t$ et $v$ non nuls, de sorte que $t^a v^b$ d\'efinit une section partout non nulle 
de $\Cal F_{a,b}$, et donc une trivialisation de ce faisceau. De m\^eme pour le type {\bf (b -)} puisque l'on peut alors supposer
$y$ et $u$ non nuls. La situation est  d'ailleurs analogue au voisinage des composantes de type (a) , mais nous
n'utiliserons pas ces derni\`eres dans la suite.
\medskip
Consid\'erons une composante du type {\bf (b +)}  et pla\c cons-nous sur le ferm\'e du voisinage
de cette composante o\`u $u=0$ : les drapeaux $(p,L)$ correspondants sont tels que $L$ passe par le 
point $p_\infty$. Un tel voisinage peut se d\'ecrire comme le quotient de l'ensemble des
$(\gth, y)$ par l'action de $\gC \cap V$ telle qu'une matrice 
 $\pmatrix 1 & \ga & \gb \\ 0 & 1 &  \bar \ga \\ 0 & 0 & 1\endpmatrix$ 
envoie $(\gth , y )$ sur $ (\exp ({2 \pi i  \over \gb_0} (\ga y + \gb)) \ \gth \ , \ y + \bar \ga ).$
\medskip
Un tel quotient s'identifie donc \`a l'espace total, que nous noterons $\tilde \Cal L_1$,
 du fibr\'e $\Cal L_1$ sur
$\Cal E$ consid\'er\'e ci-dessus. Si $\gw \in   H^1( \Gamma \setminus{ \Omega} \ , \ {\Cal F}_{a, b})$ est une
classe de cohomologie, on peut la restreidre \`a ce quotient -- ou plus pr\'ecis\'ement  \`a un 
voisinage \'epoint\'e de la section nulle dans  $\tilde \Cal L_1$ . Ensuite on  d\'efinit un
'd\'eveloppement de Fourier' comme suit.
\medskip
Notons $\Cal U \subset  \tilde \Cal L_1$ le voisinage consid\'er\'e, et $\Cal U'$ le compl\'ementaire de la section nulle
$\Cal E$. On peut prendre $\Cal U$ assez petit de sorte que le faisceau ${\Cal F}_{a,b}$  soit trivial
sur $\Cal U'$ ;  nous le trivialisons comme expliqu\'e ci-dessus.
On peut ainsi voir la restriction de $\gw$ comme un \'el\'ement de 
$$H^1({ \Cal U'}\ , \  \Cal O) = H^1({ \Cal U}\ , \ j_* (\Cal O))$$
o\`u $j$ d\'esigne l'inclusion de $\Cal U'$ dans $\Cal U$ (noter que les $R^i j_*$ sup\'erieurs sont nuls.)
La restriction \`a $\Cal E$ nous fournit un \'el\'ement de $H^1({ \Cal E}\ , \ j_* (\Cal O)_{\vert \Cal E})$.
\medskip
Ensuite, \'etant donn\'ee une section locale de $ j_* (\Cal O)_{\vert \Cal E}$, on peut \'ecrire son 
d\'eve\-loppement de Laurent sur les fibres de la projection sur $\Cal E$ (si l'on veut, en utilisant les notations 
de ci-dessus, c'est le d\'eveloppement en les puissances de $\gth$) ; les
coefficients de ce  d\'eveloppement sont des
sections locales des puissances du fibr\'e dual. On obtient donc ainsi pour chaque entier
$r \in \Bbb Z$ un morphisme de   $ j_* (\Cal O)_{\vert \Cal E}$ dans $\Cal L_{-r}$, d'o\`u une
application de $H^1({ \Cal E}\ , \ j_* (\Cal O)_{\vert \Cal E})$ dans 
$H^1({ \Cal E}\ , \Cal L_{-r})$. Au bout du compte cette construction associe \`a une classe de cohomologie
automorphe $\gw$ des 'coefficients' $\gw_r \in H^1({ \Cal E}\ , \Cal L_{-r})$ (n\'ecessairement nuls si 
$r < 0$).
\medskip
On d\'efinit de fa\c con compl\`etement analogue les coefficients de \break $\gw' \in 
H^1( \Gamma \setminus{ \Omega} \ , \ {\Cal
F}_{a, b}) $ au voisinage d'une composante  fronti\`ere de type {\bf (b -)}: ce sont des \'el\'ements
 $\gw'_r \in H^1({ \Cal E'}\ , \Cal L'_{-r})$.

\vskip1cm

\head {\bf 4. Des transformations cohomologiques au niveau des courbes elliptiques} \endhead
\bigskip
Nous allons maintenant d\'efinir des transformations cohomologiques entre les courbes elliptiques
$\Cal E $ et $\Cal E'$. Plus pr\'ecis\'ement, nous allons d\'efinir des applications lin\'eaires   (en fait bijectives)
$\eta' : H^0 (\Cal E'  , \Cal L'_r) \rightarrow H^1 (\Cal E  , \Cal L_{-r})$ et 
$\eta : H^0 (\Cal E  , \Cal L_r) \rightarrow H^1 (\Cal E'  , \Cal L'_{-r})$. Leur construction est 
analogue  \`a celle des transformations $\Cal P$ et $\Cal P'$, que nous avions rappel\'ee au d\'ebut de 
cet article.

\bigskip \noindent
{\bf (4.1)} Commen\c cons par introduire l'espace
$$\bold W = \Cal R \setminus \Bbb C \times \Bbb C $$
 o\`u $\Cal
R$ est le r\'eseau qui d\'efinit les courbes elliptiques $\Cal E$ et $\Cal E'$, op\'erant par 
$\ga \rightarrow \bar \ga$ sur le premier facteur $\Bbb C$ et par $\ga \rightarrow - \ga$ son sur
le second. Cet espace $\bold W$ se projette sur $\Cal E$ (premier facteur) et sur $\Cal E'$ (second facteur).
Nous noterons $\pi_1$ et $\pi_2$ ces deux projections.

\proclaim{Proposition 2}
 L'espace $\bold W$ est de Stein. 
\endproclaim
\demo{Preuve}  On construit comme suit des fonctions holomorphes
sur $\bold W$ : Partant de deux fonctions th\'eta $f$ et $f'$ relatives \`a $\Cal E$
et $\Cal E'$ respectivement :

$$\aligned  f (z + \bar \ga) = f (z) \ \exp ({2 \pi i r \over \gb_0} (\ga z + \gb)) , \\  \  {\text et}  \hskip 1cm
f' (z' -  \ga) = f' (z') \ \exp ({2 \pi i r \over \gb_0} (- \bar \ga z' + \gb))  , \endaligned$$
on voit alors que la fonction 
$g(z , z') = f (z) f' (z') \exp ({2 \pi i r \over \gb_0} z z')$ est invariante 
sous $\Cal R$ et d\'efinit donc une fonction sur $\bold W$. De telles fonctions s\'eparent
les points de $\bold W$ (cela r\'esulte du fait que pour $r$ assez grand les fonctions th\'eta
produisent des plongements projectifs.)
\medskip
Il reste \`a voir que $\bold W$ est holomorphiquement convexe. Donnons-nous une suite sans valeur
d'adh\'erence d'\'el\'ements de $\bold W$. Quitte \`a extraire une sous-suite, on peut supposer
des repr\'esentants $(z_n ,  z'_n)$ choisis tels que $z_n$ converge vers $z 
\in \Bbb C$. Puis (extrayant peut-\^etre encore une suite) que $z'_n = z''_n - \ga_n$ avec
$\ga_n \in \Cal R$ et $z''_n$ convergeant vers $z''$.

Il reste alors \`a choisir des fonctions th\'eta $f$  (resp. $f'$) comme ci-dessus et ne s'annulant
pas en $z$ (resp. $z''$). Alors
$$g(z_n , z'_n) = f (z_n) f' (z''_n) \exp ({2 \pi i r \over \gb_0} z_n 
z''_ n)
\exp ({2 \pi i r \over \gb_0}(- \ga_n  z_n - \bar \ga_n  z''_n  + \bar \gb_n))  \ ,$$
et on voit que cette expression n'est pas born\'ee : en effet la suite $\ga_n$ n'est pas born\'ee
et $\Re \gb_n = {1 \over 2} \vert \ga_n \vert^2$,  d'o\`u il r\'esulte que la derni\`ere exponentielle 
de l'expression ci-desssus n'est pas born\'ee (rappelons que $2 \pi i r \over \gb_0$ est un r\'eel $>0$).

\enddemo
\bigskip \noindent
{\bf (4.2)}  Nous allons maintenant d\'efinir les transformations $\eta$ et $\eta'$ en  termes 
de
la th\'eorie de Eastwood -- Gindikin -- Wong (cf. \S1).  Partons d'un \'el\'ement \hfill \break
$\phi' \in  H^0 (\Cal E'  , \Cal L'_r)$, autrement dit d'une fonction th\'eta comme plus haut.
\linebreak La fonction $\displaystyle h(z,z') =  \phi'(z') \exp ({2 \pi i r \over \gb_0} zz')$ v\'erifie 
alors l'\'equation fonctionnelle :
$$h( z+ \bar \ga , z' - \ga) = h(z,z') \exp (-{2 \pi i r \over \gb_0} (\ga z + \gb)) , $$
 c'est donc une section sur $\bold W$ de l'image r\'eciproque par $\pi_1$
du faisceau $\Cal L_{-r}$. L'ex\-pres\-sion  :
$$\eta' (\phi') = \phi'(z') \exp ({2 \pi i r \over \gb_0} zz') dz'$$
d\'efinit alors une forme diff\'erentielle relative (\`a la projection $\pi_1$)
\`a valeurs dans $\pi_1^* \Cal L_{-r}$ et donc l'\'el\'ement cherch\'e de $H^1 (\Cal E \ , \ \Cal L_{-r})$.
\medskip
On d\'efinit $\eta$ de mani\`ere  similaire : 
$$\eta (\phi) = \phi (z) \exp ({2 \pi i r \over \gb_0} zz') dz .$$

\bigskip \noindent
{\bf (4.3)}
Exprimons ces transformations en termes  de cohomologie de Dolbeault. Ceci est analogue 
\`a ce que nous avions rappel\'e plus haut (remarque \`a la fin du \S1) concernant les transformations 
$\Cal P$ et $\Cal P'$. Le moyen d'effectuer cette  traduction 
est expliqu\'e dans \cite{EGW} : pour obtenir un repr\'esentant de la classe de 
cohomologie de Dolbeault associ\'ee \`a $\eta' (\phi')$, on doit commencer par \'etendre cette derni\`ere 
en une forme diff\'erentielle absolue sur $\bold W$; puis prendre l'image r\'eciproque de 
cette extension par $s$, une section $C^\infty$ de la projection de $\bold W$ sur $\Cal E$. Enfin,
on prend la partie de type $(0,1)$ de cette image r\'eciproque.

\medskip
Dans notre cas une
section est donn\'ee par $s(z) = (z , -\bar z)$;  si l'on  prend 
la partie de type $(0,1)$ du pull-back par $s$ de la forme diff\'erentielle 
$\phi' (z') \exp ( {2 \pi  i r
\over \gb_0}  z z') \  dz'$, on obtient l'expression de notre \'el\'ement comme une
classe de cohomologie de Dolbeault:
$$\eta'_{\Cal D} (\phi') = - \phi' (- \bar z) \ \exp ( - {2 \pi  i r
\over \gb_0}  z \bar z) \  d \bar z \ .$$

De fa\c con analogue $\eta (\phi)$ s'exprime en cohomologie de Dolbeault :
$$\eta_{\Cal D} (\phi) = - \phi (- \bar z') \ \exp ( - {2 \pi  i r
\over \gb_0}  z' \bar z') \  d \bar z' \ .$$

\bigskip \noindent
{\bf (4.4)}
La formule 
$$\langle \phi'_1 , \phi'_2 \rangle = \int_{\Cal E'} \phi'_1 (z') \overline  {\phi'_2 (z')}
\ \exp (- {2 \pi i r \over \gb_0} \vert z' \vert^2 ) \  \ \vert d z' \wedge \overline { d z'} \vert $$
d\'efinit un produit scalaire hermitien sur  $H^0 (\Cal E' , \Cal L'_ r)$, associ\'e \`a 
une m\'etrique sur le fibr\'e $\Cal L'_r$. Cette m\'etrique est d'ailleurs essentiellement
canonique (i.e. \`a une constante de normalisation pr\`es) :  elle est caract\'eris\'ee,
\`a multiplication pr\`es par un scalaire, par le fait d'\^etre  hermitienne et \`a courbure parall\`ele.
 D'autre part, on a
une bijection  antilin\'eaire  de $H^0 (\Cal E' , \Cal L'_ r)$ sur $H^0 (\Cal E , \Cal L_ r)$ qui associe \`a
$\phi'$ la section $\check \phi'$ d\'efinie par $\check \phi'(z) = - \overline {\phi' (- \bar z)}$.
On a finalement la  dualit\'e de Serre $(Ê\ , \ )$ entre $H^0 (\Cal E , \Cal L_ r)$ et 
$H^1 (\Cal E , \Cal L_ {-r})$, donn\'ee en cohomologie de Dolbeault par l'int\'egrale sur $\Cal E$
du produit.
Il r\'esulte des formules ci-dessus que pour $\phi'_1$ et $\phi'_2 $ $\in  H^0 (\Cal E' , \Cal L'_
r)$ on a
$$\langle \phi'_1 , \phi'_2 \rangle = (\eta'_{\Cal D} (\phi'_1) , \check \phi'_2 ) .$$
ce qui exprime l'application $\eta'$ comme une certaine compos\'ee entre dualit\'e de
Serre, dualit\'e $\langle \ ,  \ Ê\rangle$, et application $\phi'\rightarrow \check \phi'$.
En particulier $\eta'$ est une bijection entre $H^0 (\Cal E' , \Cal L'_ r)$ et
$H^1 (\Cal E , \Cal L_ {- r})$. De m\^eme bien s\^ur pour $\eta$.

\bigskip \noindent
{\bf (4.5)}  Les courbes $\Cal E$ , $\Cal E'$ ainsi que les faisceaux $\Cal L$ et $\Cal L'$ sont
d\'efinis  sur l'extension ab\'elienne maximale $F^{ab} \subset  \Bbb C$ de $F$. 
Notre objectif est  de prouver qu'il en est de m\^eme des transformations $\eta$ et 
$\eta'$.

\proclaim {Th\'eor\`eme 1 }  A multiplication pr\`es par un scalaire ind\'ependant de $r$,
les  transformations $\eta$ et 
$\eta'$ sont d\'efinies sur $F^{ab}$ (autrement dit transforment sections ration\-nelles sur
ce corps en classes de cohomologie rationnelles).
\endproclaim
\demo{Preuve} D'apr\`es ce qu'on a dit plus haut (4.4), il suffit de voir  que 
 le produit scalaire hermitien sur  $H^0 (\Cal E' ,
\Cal L'_ r)$ d\'efini par la formule ci-dessus, de m\^eme que son analogue sur 
   $H^0 (\Cal E , \Cal L_ r)$:

$$\langle \phi_1 , \phi_2 \rangle = \int_{\Cal E} \phi_1 (z) \overline  {\phi_2 (z)}
\ \exp (- {2 \pi i r \over \gb} \vert z \vert^2 ) \ \  \vert  d z \wedge \overline { d z} \vert $$
est {\it rationnel} sur $F^{ab}$ \`a une constante multiplicative $c$ pr\`es,
c'est-\`a-dire que si $\phi_1$ et $\phi_2$ sont d\'efinis sur $F^{ab}$ alors il en 
est de m\^eme de $c\langle \phi_1 , \phi_2 \rangle$. 
\medskip
On utilise dans ce but l'article \cite{Sh1} qui explicite une base de telles sections
rationnelles sur $F^{ab}$ en termes de fonctions th\'eta classiques. On peut alors 
explicitement calculer les produits scalaires de ces \'el\'ements, dont on trouve qu'ils
constituent une base orthogonale. Un tel calcul se trouve d'ailleurs dans Siegel \cite{Si}.
Tout cela vaut d'ailleurs plus g\'en\'eralement pour des vari\'et\'es ab\'eliennes CM.
\medskip
Notons $(\gw_1, \gw_2)$ une base du r\'eseau $\Cal R$ choisie de telle sorte que 
$z = \gw_1 / \gw_2$ soit de partie imaginaire $>0$. Posons $\gl = { i r \over \gb_0}$.
Alors les fonctions th\'eta que nous consid\'erons sont celles \'etudi\'ees dans loc.cit.
avec $H (x , y) = 2 \gl \bar x y$ (\`a ne pas confondre bien s\^ur avec la forme hermitienne
du \S1 qui d\'efinit $G$) et $\Psi (a) = \exp (2 \pi r i {\Im b \over \Im \gb_0})$
(une racine de l'unit\'e). Les notations de Shimura s'\'ecrivent dans notre cas particulier :
 $\gW = (\gw_1, \gw_2)$,  $\epsilon = 1$ et enfin
$$\mu = E ((\gw_2, \gw_1) = 2 \gl \Im ( \bar \gw_2
\gw_1) = 2 \gl \vert \gw_2 \vert^2 \Im z. $$ 
\medskip
Les fonctions th\'eta sont donn\'ees par 
$$\theta (u,z ; \gr,s) = \sum_{n \in \Bbb Z} \exp (\pi i (n+\gr)^2 z) \exp (2 \pi i (n+\gr) (u+s)) ,$$
o\`u $u\in \Bbb C$ et o\`u $\gr$ et $s$ sont des param\`etres r\'eels (en fait dans $\Bbb Q$) qui doivent \^etre choisis de
telle sorte que $\psi (\gw_1 \ga_1 + \gw_2 \ga_2) = \exp (\pi i \ga_1 \ga_2 +  \mu \gr \ga_2
- s \ga_1).$ Le rationnel $s$ est d\'efini   \`a un entier pr\`es, et $\gr$ \`a un \'el\'ement pr\`es
de
$\gm^{-1} \Bbb Z$.
\medskip
On pose ensuite $\phi (u,z ; \gr , s) =  \exp ({\pi \over 2 \Im z} u^2) \ \theta (u,z ; \gr , s)$
et on d\'efinit enfin des sections de $\Cal L_r$ par la formule :
$$f_{\gr,s} (u) = \phi ( \gm \gw_2^{-1} u , \gm \gw_2^{-1} \gw_1 ; \gr , s) 
= \exp ({\pi \over 2 \Im z} \mu \gw_2^{-2} u^2) \ \theta (\gm \gw_2^{-1} u , \gm z ; \gr , s);$$
Soit encore, avec $\zeta = \gw_2 / \vert \gw_2 \vert$ :
$$f_{\gr , s} (u) = 
 \exp (\pi \gl  \zeta^{-2} u^2) \ \theta (\gm \gw_2^{-1} u , \gm z ; \gr,s) .$$
\medskip
Un des r\'esultats de l'article pr\'ecit\'e est le suivant: Pour $\gr, s$ fix\'es et $j$ variant
dans un syst\`eme de repr\'esentants des classes de $\gm^{-1} \Bbb Z / \Bbb Z$, les fonctions
$f_{\gr + j \  , \  s} (u)$ constituent une base de $H^0 (\Cal E , \Cal L_ r)$ sur $\Bbb C$. D'autre part 
on obtient une base de sections rationnelles sur $F^{ab}$ en prenant les fonctions 
$f'_{\gr + j  \ ,  \ s} (u) =  h(z)^{-1} f_{\gr + j  \  , \ s} (u)$ pour $h$ une quelconque forme modulaire de poids
demi-entier alg\'ebrique sur $\Bbb Q^{ab}$.
\medskip
On voit que les produits scalaires $\langle f_{\gr_1 , s} f_{\gr_2 , s} \rangle $ sont donn\'es
par les int\'egrales :
$$
\int_{\Cal E} \exp(-4 \pi \gl (\Im (\zeta^{-1} u))^2  \theta (\gm \gw_2^{-1} u , \gm z ; \gr_1,s) 
\overline {\theta (\gm \gw_2^{-1} u , \gm z ; \gr_2,s)} \ \vert \ d u
\wedge \overline { d u} \vert $$
\medskip
On calcule ces derni\`eres en rempla\c cant les fonctions th\'eta par la somme qui les d\'efinit
et en int\'egrant sur un rectangle fondamental (de sommets $0 , \gw_1, \gw_2,
\gw_1 + \gw_2$) chacun des produits
qui apparaissent. Autrement dit on doit int\'egrer
$$ \exp \big(2 \pi i (n_1 +\gr _1) (\gm \gw_2^{-1} u+s) - 2 \pi i (n_2 +\gr _2) (\gm \bar \gw_2^{-1}
\bar u + s) - 4 \pi \gl \vert \gw_2 \vert^2 (\Im (\gw_2^{-1} u))^2\big) = $$
$$\exp  \big(2 \pi i (n_1 + \gr_1 - n_2 - \gr_2)\gm  \Re (\gw_2^{-1} u) - 2 \pi  (n_1 + \gr_1 + n_2 +
\gr_2) \gm \Im (\gw_2^{-1} u)  - 4 \pi \gl \vert \gw_2 \vert^2 (\Im (\gw_2^{-1} u))^2\big)$$
\medskip
On peut commencer cette int\'egration avec $\Im (\gw_2^{-1} u)$ constant et la partie r\'eelle 
variant de $0$ \`a $1$ ;  les int\'egrales ci-dessus sont alors nulles sauf si 
$n_1 + \gr_1 - n_2 - \gr_2 = 0$. Ceci n'est possible que pour $ \gr_1 \equiv \gr_2 \mod \Bbb Z$, et donc
dans notre cas on voit que les fonctions que nous consid\'erons sont orthogonales deux \`a deux.
\medskip
Il reste \`a calculer   $\langle f_{\gr , s} f_{\gr , s} \rangle $.  C'est la somme de la s\'erie 
dont le terme de rang $n$ est   le produit par $\exp (- 2 \pi  (n+\gr)^2 \gm \Im z)$ de l'nt\'egrale
sur le rectangle fondamental de $\exp  \big( - 4 \pi  (n + \gr)\Im (\gw_2^{-1} u)  - 4 \pi \gl
\vert \gw_2 \vert^2 (\Im (\gw_2^{-1} u))^2\big)$. On doit donc calculer la somme des int\'egrales
 des fonctions :
$$\exp \big (- 2 \pi  (n+\gr)^2 \gm \Im z  - 4 \pi  (n + \gr) \gm \Im (\gw_2^{-1} u)  - 4 \pi \gl
\vert \gw_2 \vert^2 (\Im (\gw_2^{-1} u))^2\big) .$$
Soit, compte tenu du fait que $\gm \Im z = 2 \gl  \vert \gw_2 \vert^2 (\Im z)^2$ , la somme des
int\'egrales des $\exp \Big( - 4 \pi \gl \vert \gw_2 \vert^2 \big( (n+\gr) \Im z + \Im (\gw_2^{-1}
u)
\big)^2
\Big)$.

Prenons pour fixer les id\'ees la mesure usuelle sur $\Bbb C$. Le changement de variable 
$v = \gw_2^{-1} u$ transforme notre rectangle en celui de sommets $0$, $1$, $z$, $1 + z$, et nous
devons donc int\'egrer sur ce dernier les fonctions
 $\vert \gw_2 \vert^2 \exp \Big( - 4 \pi \gl \vert \gw_2 \vert^2 \big( (n+\gr)
\Im z + \Im (v) \big)^2 \Big)$. Comme cette fonction ne d\'epend que de $y = \Im u$, on doit donc
calculer la somme des $\vert \gw_2 \vert^2 \int_0^{\Im z}\exp \Big(- 4 \pi \gl \vert \gw_2
\vert^2\big( (n+\gr)
\Im z + y \big)^2 \Big) dy $ . On  voit alors finalement que l'on a:
 $$\langle f_{\gr , s} f_{\gr , s} \rangle = \vert \gw_2 \vert^2 \int_{-\infty}^{+ \infty}
\exp(- 4 \pi \gl \vert \gw_2 \vert^2 y^2) dy$$
\medskip \noindent 
On obtient donc le r\'esultat suivant :
$$ \langle f_{\gr , s} f_{\gr , s} \rangle = {\vert \gw_2 \vert \over 2 \sqrt \gl} $$
et donc 
$$ \langle f'_{\gr , s} f'_{\gr , s} \rangle = {\vert \gw_2 \vert \over 2  \vert  h(z)^2 \vert   \sqrt \gl}  \ ;$$
ces produits scalaires sont donc bien rationnels sur $F^{ab}$ \`a la constante pr\`es \break
$c = \vert {h(z)^2 \over \gw_2 } \vert$.
\enddemo

\vskip1cm
\head {\bf 5. Relations entre les transformations $\Cal P , 
\Cal P', \eta , \eta'$}
\endhead
\bigskip
\proclaim{Th\'eor\`eme  2}
 Soit $f'$ une forme modulaire parabolique de poids $k \geq 2$ sur 
$\Gamma \backslash   \overline \Delta$ et soit 
$\gw' =\Cal P' (f') \in  H^1( \Gamma \setminus{ \Omega} \ , \ 
{\Cal F}_{-k +1, k-2})$ la
classe de cohomologie automorphe correspondante. Alors les coefficients du
d\'eveloppement  de $\gw'$ au voisinage de la composante fronti\`ere du type 
{\bf ( b -) } sont nuls, tandis qu'en celle du type {\bf ( b +) }
ces coefficients sont donn\'es en fonction des coefficients $g'_r$
de $f'$ par : $\gw'_r = - \eta'(g'_r)$ pour $r \geq 1$, et  $\gw'_r = 0$
pour $r \leq 0$ .

De m\^eme pour $f$ une telle forme sur $\Gamma \backslash   
\Delta$, les coefficients de $\gw = \Cal P (f)$ sont donn\'es 
par $\gw_r = - \eta(g_r)$ ($r \geq 1$) et $0$ si $r \leq 0$ au voisinage de la composante de type 
{\bf (b -)} ;  ils sont  nuls au voisinage de celle de type {\bf (b +)}.
\endproclaim
\demo{Preuve}
\medskip
(i) Rappelons en quelques mots la construction 
de \cite{EGW};  Soit  $\pi : X \rightarrow M$ une fibration holomorphe
entre vari\'et\'es complexes, \`a fibres contractiles et dont l'espace
total est de Stein. Soit
$\Cal F$ un fibr\'e vectoriel sur $M$. Alors on a  un isomorphisme entre 
la cohomologie $H^* (M , {\Cal O} (\Cal F))$ et la cohomologie du
complexe 
$\gC (X , \gW^*_\pi (\Cal F))$, o\`u  $\gW^*_\pi (\Cal F)$ d\'esigne le faisceau
des formes diff\'erentielles $\pi$--relatives \`a valeurs dans $\Cal F$. Cet isomorphisme
s'obtient ainsi : on consid\`ere le
faisceau
$\pi^* {\Cal O} (\Cal F)$ image r\'eciproque de ${\Cal O} (\Cal F)$.
Parce que les fibres de $\pi$ sont contractiles, on a : $H^* ((M , {\Cal
O} (\Cal F))
\simeq H^* ( X, \pi^* {\Cal O} (\Cal F)) $.  D'autre part le complexe
de de Rham relatif est une r\'esolution du faisceau $\pi^* {\Cal O}
(\Cal F)$, et cette r\'esolution est 
acyclique parce que l'espace total est de Stein.
D'o\`u un isomorphisme entre  $H^* ( X, \pi^* {\Cal O} (\Cal F)) $ et 
$H^* ( \gC (X , \gW^*_\pi (\Cal F)))$.
\medskip
On voit facilement que cette construction est fonctorielle au
sens suivant : pour $M' \rightarrow M$ et $\eta$ une classe de cohomologie
 associ\'ee
\`a une forme diff\'erentielle $\gw$, sa restriction \`a $M'$ est
donn\'ee par la restriction de $\gw$ au produit fibr\'e 
$X' = M' \times_M X$ pourvu que ce dernier soit de Stein. C'est le cas
en particulier pour $M'$ un sous-espace ferm\'e de $M$. 
\medskip
De m\^eme, si une classe de cohomologie est associ\'ee \`a $\gw$ sur
$X$ et si $X' \subset X$ est une sous-vari\'et\'e, encore de Stein, et
telle que la restriction de $\pi$  soit encore une fibration sur $M$ \`a 
fibres contractiles, alors la m\^eme classe est d\'efinie par la restriction de
$\gw$ \`a $X'$.

\medskip

(ii) Nous appliquerons ces consid\'erations  \`a $\Cal P'(f')$, qui est d\'efinie par la
forme diff\'erentielle $  f' ( z \wedge \xi) \ l(\xi)^k \ \omega_J
\  $ sur l'espace  $\Gamma \setminus \bold U$ (cf.  \S1 ; on tient compte
ici du fait que $a=0$ ). En fait nous commencerons par prendre  le pull-back de cette 
forme sur l'espace quotient  $\Gamma \cap V \setminus \bold U$ : ceci
d\'efinit  le pull-back sur $\Gamma \cap V \setminus \gW$ de notre classe 
de cohomologie ; noter que  $\Gamma \cap V \setminus \bold U$ est de Stein : 
cela r\'esulte de ce qu'il est un rev\^etement de l'espace de Stein
$\Gamma \setminus \bold U$ (on peut aussi  d\'emontrer cela de fa\c con plus explicite).

\medskip
Notons $\bold U^-$ sous-ensemble ferm\'e de $ \bold U$ constitu\'e des couples de
drapeaux
$(z,l;\xi
\ga)$ dans
$\bold U$ tels que $z \in D_\infty $ et que $\ga$ passe par
$p_\infty$.  Le quotient $\gC \cap V \setminus \bold U^-$, ferm\'e dans 
$\gC \cap V \setminus \bold U$, est encore de Stein. Il se projette sur le sous ensemble
ferm\'e de 
$\gC \cap V \setminus \gW$, que nous noterons $S^-$ et qui est  constitu\'e des drapeaux
qui v\'erifient la condition  suppl\'ementaire $z \in D_\infty $.  Le d\'eveloppement 
en une composante de type {\bf ( b -) } se calcule \`a partir de 
la restriction \`a $S^-$ de la classe de cohomologie $\gw$.

On voit que les fibres de la  projection de $\gC \cap V \setminus \bold
U^-$ sur $S^-$ sur sont encore
contractiles. Par suite la restriction de notre classe de cohomologie est
donn\'ee par la restriction de la forme diff\'erentielle correspondante \`a 
$\gC \cap V \setminus \bold U^-$.

De fa\c con duale on consid\`ere $\bold U^+$, form\'e des drapeaux tels
que $l$ passe par $p_\infty$ et que $\xi \in D_\infty$. C'est un rev\^etement
de sous-espace $S^+ \subset \gC \cap V \setminus \gW$ constitu\'e des drapeaux
$(z,l)$ v\'erifiant la condition suppl\'ementaire que $l $ passe par $p_\infty$.
Le d\'eveloppement en la composante de type  {\bf ( b +) } se lit sur la restriction de 
la classe de cohomologie \`a $S^+$, et cette restriction est associ\'ee \`a la 
restriction \`a  $\bold U^+$ de la forme diff\'erentielle.

\medskip
(iii) Commen\c cons par le cas  d'une composante 
de type {\bf ( b -) }. On va montrer alors que la restriction de  $\Cal P'(f')$
\`a $\bold U^-$ est relativement exacte.  Pour cela, on  remarque 
que c'est vrai sur chaque fibre de la projection sur $S^-$ : $\Cal P'(f')$
y est la
diff\'erentielle d'une section holomorphe bien d\'efinie \`a une constante 
pr\`es. Ces int\'egrales sur les fibres ne d\'ependent que de $J$, la droite 
joignant $x$ \`a $\xi$, et nous allons  les normaliser de sorte que leur limite soit
nulle lorsque cette droite tend vers $D_\infty$.  Une fois cette normalisation
accomplie, on obtient une section  holomorphe bien d\'efinie sur  $\bold U^-$
du faisceau $\Cal F_{-k + 1  \ k-2}$ dont  $\Cal P'(f')$ est la diff\'erentielle
relative.

V\'erifions qu'une telle normalisation est possible. Soient 
$z = \pmatrix x_0 \\ 1 \\ 0 \endpmatrix$ et  \break $l = \pmatrix 1 & -x_0 & \gt_0 \endpmatrix $
un point fix\'e de $S^-$ ; un point de la fibre au-dessus de cet \'el\'ement 
correspond \`a la donn\'ee suppl\'ementaire de 
$\xi = \pmatrix x \\ y \\ 1 \endpmatrix$ et
 $\ga = \pmatrix 0 & 1 & -y \endpmatrix$ tels que le point d'intersection $I$ de 
$l$ et $\ga$ soit int\'erieur \`a $\gD$ et que la droite $J$ joignant $z$ et $\xi$
soit ext\'erieure \`a $\overline \gD$. Des coordonn\'ees de $J$ sont donn\'ees par
$J =   \pmatrix 1 & -x_0 & w \endpmatrix $ avec $w = x_0 y - x $ ( 
$2 \Re w > \vert x_0 \vert^2$). La
diff\'erentielle
$\gw_J$ est la diff\'erentielle d'une coordonn\'ee $y(J)$ bien d\'efinie \`a une 
constante pr\`es par : $y(J) = \det_z (J,J_0 ) \det_z^{-1} (J, l)$, o\`u les
d\'eterminants sont pris dans le plan constitu\'e des formes nulles sur $z$
et o\`u $J_0$ repr\'esente  le choix d'un point-base (ne d\'ependant que de $(z,l)$),
 normalis\'e de
telle sorte que
 $\det (J_0, l , * ) = * (z)$. 
\medskip
 Ici on peut prendre $J_0 = \pmatrix 0 & 0 & 1 \endpmatrix $, ce qui donne
$y(J) = (\gt_0 - x_0 y + x )^{-1} = (\gt_0 - w)^{-1}$.  Parce que 
$l(\xi) = x - x_0 y + x = \gt_0 - w$, la restriction \`a la fibre au dessus de $(z,l)$ de notre
forme diff\'erentielle  est donn\'ee par :
 $$\Cal P'(f') = f' (J) (\gt_0 - w )^k d((\gt_0 - w)^{-1} = f'(1, \  -x_0,\  w)
 \  \ (\gt_0 - w)^{k-2} dw \ .$$ 

L'expression $f' ( 1, \ -x_0 , \ w)$ est une fonction $\phi'(q)$ de 
$q = \exp (-{2 \pi i \over \gb_0} w)$, holomorphe au voisinage de $0$ et nulle en 
$0$. Notre forme diff\'erentielle est alors le produit de ${\phi'(q)\over q } dq$ par un
polyn\^ome en $\log (q)$ et on voit qu'une primitive d'une telle forme admet une limite
en $q = 0$ (une d\'etermination de $\log$ \'etant choisie avec une partie imaginaire
born\'ee); autrement dit, toute primitive admet une  limite lorsque $\Re w$ tend vers
$+\infty$, la  partie imaginaire restant born\'ee.

\medskip
(iv)
 Traitons maintenant  le cas  d'une composante 
de type {\bf ( b +) }. Un point de $\bold U^+$ \'etant d\'efini par un couple de drapeaux
  $ (z,l) =(\pmatrix x  \\ y 
 \\  1 \endpmatrix , \pmatrix 0 & 1 & -y \endpmatrix) $ et $ (\xi ,\ga) = (\pmatrix -v
\\ 1 \\ 0
\endpmatrix ,
 \pmatrix 1 & v & w \endpmatrix)$, on voit que la  droite $J$ est donn\'ee par
 $J =
\pmatrix 1 & v & -x - vy
\endpmatrix$. La classe de cohomologie 
${\Cal P}'(f')$ est associ\'ee \`a la  forme diff\'e\-rentielle relative  
$f'(J) \  \gw_J$ (noter en effet que $l(\xi)$ est identiquement $1$).

On calcule comme ci-dessus (la restriction de) $\gw_J$. 
On peut prendre $J_0 = \pmatrix 1 &v_0 & -x - v_0 y \endpmatrix$ (avec $v_0$ fix\'e). On
trouve alors que
$y(J) = v_0 - v$ et donc  $\gw_J = - dv$.

 Finalement,  $-{\Cal P}' (f')$ est donn\'ee par:

$$-{\Cal P}' (f') (x,y;v,w)  = f' (v, -x-vy) \ dv = \sum_{r \in \Bbb N^*} g'_r (v) exp ( {2 \pi  i r
\over \gb_0} (x + vy)) \ dv$$

Ou, si l'on pr\'ef\`ere, en termes $\gth = \exp ({2 \pi i  \over \gb_0}
x)$ :

$$-{\Cal P}' (f') (x,y;v,w)   = \sum_{r \in \Bbb N^*} g'_r (v) exp ( {2 \pi  i r
\over \gb_0}  vy) \  \gth^r \ dv \ .$$

\medskip
Cette derni\`ere expression ne d\'epend que de $v$, $y$ et $\gth$. L'action 
d'un \'el\'ement $\pmatrix 1 & \ga & \gb \\ 0 & 1 &  \bar \ga \\ 0 & 0 & 1\endpmatrix$
de $\gC \cap V$ transforme $y$ en $y + \bar \ga$, $v$ en $v - \ga$, $\gth $ en
$\exp {2 \pi i \over \gb_0} (\ga y + \gb)$, et laisse invariante l'expression de 
${\Cal P}' (f')$.
\medskip

Notons $\bold W \tilde \Cal L_1$ l'espace total du fibr\'e image r\'eciproque de $\Cal
L_1$ par la projection de $\bold W$ sur $\Cal E$. Cet espace ce projette sur l'espace
 total du fibr\'e $\Cal L_1$, que nous avions not\'e $\tilde \Cal L_1$. Notons 
\'egalement  $\bold W \tilde \Cal L_1 ^*$ et $\tilde \Cal L_1^*$ les compl\'ementaires 
des sections nulles. On voit que l'expression pr\'ec\'edente d\'efinit une forme
diff\'erentielle sur $\bold W \tilde \Cal L_1 ^*$ (relativement \`a sa projection 
sur $\tilde \Cal L_1^*$).
\medskip
Remarquons de plus que l'espace  $\bold W \tilde \Cal L_1$ est de Stein : on peut
utiliser le fait g\'en\'eral qu'un fibr\'e vectoriel sur un espace de Stein est encore
de Stein, ou bien aussi le prouver directement dans ce cas particulier. Il en est donc de m\^eme
de 
 $\bold W \tilde \Cal L_1^*$, et il en r\'esulte que la formule ci-dessus d\'efinit
une classe de cohomologie sur $\tilde \Cal L_1^*$. Cette derni\`ere n'est autre, en vertu
des consid\'erations (i) ci-dessus, que la classe de cohomologie qui nous a servi en (3.2) 
\`a d\'efinir le d\'eveloppement de Fourier.

\medskip 
On remarque enfin que la formule ci-dessus a un sens sur l'espace
$\bold W \tilde \Cal L_1$ tout entier : notre classe de cohomologie 
se prolonge donc \`a $\tilde \Cal L_1$. Dans ces conditions il est clair que la
construction de (3.2) revient \`a se restreindre \`a la section nulle puis 
\`a prendre les coefficients de Fourier. Il en r\'esulte 
que la $r$--i\'eme composante de Fourier de $-{\Cal P}' (f')$ est  associ\'ee au
cocycle qui envoie
$(y,v)$ sur $g'_r (v) exp ( {2 \pi  i r
\over \gb_0}  vy) \ dv$. Cela prouve le th\'eor\`eme pour la transformation $\Cal P$.
\bigskip
La d\'emonstration dans le cas de la transformation $\Cal P$ est duale de 
la pr\'ec\'edente.
\enddemo
\bigskip
Compte tenu des propri\'et\'es de rationalit\'e pour les transformations $\eta$
et $\eta '$ prou\-v\'ees au paragraphe pr\'ec\'edent, et du fait que les formes modulaires 
de Picard sont engendr\'ees par celles rationnelles sur $F^{ab}$, le th\'eor\`eme
pr\'ec\'edent admet le corollaire suivant : 
\medskip
\proclaim{\bf  Corollaire 1}  L'image de la transformation $\Cal P'$ est engendr\'ee par des classes
de cohomologie rationnelles sur $F^{ab}$ (au sens que leurs coefficients sur les composantes du type
(b+) le sont). R\'esultat analogue pour l'image de $\Cal P$.
\endproclaim
On remarquera d'ailleurs que cette image est l'ensemble des classes paraboliques, que l'on peut d\'efinir
de fa\c con analytique (remarque finale du \S1) \`a d\'efaut de disposer d'une d\'efinition g\'eom\'etrique.

\bigskip \noindent
{\bf Commentaire final : } Pour $\gw' =\Cal P' (f')$ comme dans l'\'enonc\'e du th\'eor\`eme, le fait que les coefficients
correspondant aux composantes de type (b-)  soient nuls traduit le fait qu'au voisinage d'une composante fronti\`ere de
ce type, la restriction de $\gw'$ \`a la partie alg\'ebrique de la compactification (cf. (3.2)) s'annule.  On aimerait pouvoir d\'efinir
des invariants plus fins, qui "mesurent" le $H^1$ d'un voisinage de la courbe fronti\`ere (rappelons  (cf. (2.5)) qu'un tel voisinage n'est
pas de Stein). 

Ce probl\`eme est intimement li\'e \`a celui de d\'efinir des coefficients pour les 2-classes de cohomologie du type de Maass; en effet,
de m\^eme que nous l'avions d\'emontr\'e dans \cite{Ca1} dans le cas compact, il devrait encore \^etre vrai que ces classes (dans le cas parabolique)
sont engendr\'ees par des cup-produits ${\Cal P}(f) \cup {\Cal P'}(f')$.

\enddocument